\documentclass[12pt]{amsart}
\usepackage{amsmath,amssymb,amscd}
\usepackage{stmaryrd}
\usepackage{amsfonts}
\usepackage{latexsym}
\usepackage{epsfig}
 \numberwithin{equation}{section}

\newcommand{\overbar}[1]{\mkern 1.5mu\overline{\mkern-1.5mu#1\mkern-1.5mu}\mkern 1.5mu}
\newcommand*\xoverline[2][0.75]{%
    \sbox{\myboxA}{$\m@th#2$}%
    \setbox\myboxB\null% Phantom box
    \ht\myboxB=\ht\myboxA%
    \dp\myboxB=\dp\myboxA%
    \wd\myboxB=#1\wd\myboxA% Scale phantom
    \sbox\myboxB{$\m@th\overline{\copy\myboxB}$}%  Overlined phantom
    \setlength\mylenA{\the\wd\myboxA}%   calc width diff
    \addtolength\mylenA{-\the\wd\myboxB}%
    \ifdim\wd\myboxB<\wd\myboxA%
       \rlap{\hskip 0.5\mylenA\usebox\myboxB}{\usebox\myboxA}%
    \else
        \hskip -0.5\mylenA\rlap{\usebox\myboxA}{\hskip 0.5\mylenA\usebox\myboxB}%
    \fi}

\newtheorem{lemma}{Lemma}[section]
\newtheorem{corollary}[lemma]{Corollary}
\newtheorem{theorem}[lemma]{Theorem}

\newtheorem{proposition}[lemma]{Proposition}

\theoremstyle{remark}

\theoremstyle{definition}
\newtheorem{definition}[lemma]{Definition}

\def\ch{\operatorname{ch}}
\def\td{\operatorname{td}}

\def\A{\mathcal{A}}

\def\w{\omega}
\def\O{\mathcal{O}}
\def\T{\Theta}

\def\E{\mathcal{E}}

\def\P{\mathbb{P}}

\def\Z{\mathbb{Z}}

\def\a{\alpha}

\def\d{\delta}

\def\l{\lambda}
\def\O{\mathcal{O}}

\def\s{\sigma}

\def\la{\langle} 
\def\ra{\rangle}

\def\H{\operatorname{H}}

\def\det{\operatorname{det}}
\def\GL{\operatorname{GL}}
\def\PSL{\operatorname{PSL}}

\def\B{\mathcal{B}}
\def\Cl{\operatorname{Cl}}
\def\H{\mathcal{H}}
\def\T{\mathcal{T}}
\def\M{\mathcal{M}}
\def\C{\mathbb{C}}

\def\E{\mathcal{E}}
\def\ker{\operatorname{ker}}
\def\Hom{\operatorname{Hom}}

\def\Ext{\operatorname{Ext}}

\def\End{\operatorname{End}}
\def\Pic{\operatorname{Pic}}
\def\Spec{\operatorname{Spec}}

\def\Br{\operatorname{Br}}
\def\IsGr{\operatorname{IsGr}}

\def\cal{\mathcal}
\def\L{\mathcal{L}}
\def\Sym{\operatorname{Sym}}
\def\PGL{\operatorname{PGL}}
\def\Gr{\operatorname{Gr}}
\def\Grass{\operatorname{Gr}}
\def\rk{\operatorname{rk}}

\def\st{\,|\,}

\def\GL{\operatorname{GL}}
\def\SL{\operatorname{SL}}

\def\G{\mathbb{G}}
\def\cEnd{\operatorname{\mathcal{E}\!\mathit{nd}\,}}
\newcommand{\suchthat}{\, \mid \,} % nice "such that"
\def\BSV{\operatorname{BSV}}
\def\BS{\operatorname{BSV}}
\def\mid{\,|\,}
\begin{document}

\author{Colin Ingalls}
\address{Department of Mathematics and Statistics, University of New Brunswick,
Fredericton, New Brunswick,
Canada.}\email{cingalls@unb.ca}
\author{Madeeha Khalid}
\address{Department of Mathematics,
St. Patrick's College Drumcondra, Dublin, Ireland}
\email{madeeha.khalid@spd.dcu.ie}
\thanks{
The first author was supported by an NSERC Discovery Grant,
and the second author was supported in part by
the Embark Initiative program by the Irish Research Council for Science Engineering and Technology and in part by St. Patrick's College Drumcondra.}

\title[A Derived Equivalence of Azuamaya Algebras on K3 Surfaces]{An explicit Derived Equivalence of Azumaya Algebras on K3 Surfaces
via Koszul Duality}

\begin{abstract}

We consider moduli spaces of Azumaya algebras on K3 surfaces and construct an example. In some cases we show a derived equivalence which corresponds to a derived equivalence between twisted sheaves. 
We prove if $A$ and $A'$ are Morita equivalent Azumaya algebras of degree $r$ then $2r$ divides $c_2(A) - c_2(A')$.  In particular 
this implies that
if $A$ is an Azumaya algebra on a K3 surface and $c_2(A)$ is within $2r$ of its minimal
bound then the moduli stack of Azumaya algebras with the same underlying gerbe, if non empty, is a proper algebraic space.
\end{abstract}

 \maketitle

\begin{section}{Introduction}
In this paper we investigate specific examples of moduli spaces of Azumaya algebras on K3 surfaces and their derived equivalences.

Often in mathematical physics one is interested in field theories with non trivial B-fields.  The B-field is a finite order element in the second cohomology.  Equivalences between physical theories with B-fields, can be interpreted in terms of derived equivalences of twisted sheaves \cite{Pantev-Sharpe}.  Twisted sheaves correspond to modules over Azumaya algebras when the Picard group is torsion free, and so examples of derived equivalences between modules over Azumaya algebras are also relevant to string theory.

%Some examples of K3 surfaces are; double covers of $\P^2$ branched along a smooth
%sextic curve, smooth quartic hypersurfaces in $\P^3$, complete intersections of a quadric and cubic hypersurface in $\P^4$, triple intersection of quadric hypersurfaces in $\P^5$.
%These are all algebraic K3s.
%In general for each $n$ there is a 19 dimensional moduli space of K3 surfaces occurring as normal surfaces of degree $2n-2$ in $\P^n$.

%Another interesting class of K3 surfaces called "Kummer surfaces" is realised by taking the
%quotient space of the canonical involution on a 2-dimensional complex torus $T$ and blowing up the 16 singular points. These Kummer surfaces may be non algebraic if the
%complex torus is non algebraic.  

%The paper is roughly divided into several parts.
Section 2 contains the preliminaries and basic definitions.
In section 3, we review a classical example of an Azumaya algebra $A$ that arises from a quadric bundle on $\P^2$.  The centre of $A$  is the K3 surface $M$ which is a double cover of $\P^2$ branched over a sextic curve.  The sextic curve parameterises the degenerate quadrics.  We calculate invariants of $A$ and prove that $A$ has no deformations.  Moreover $A$ has index 2 and degree 4 so is clearly not a division algebra at the generic point.  By a reduction technique we construct a family of Azumaya algebras that are Morita equivalent to $A$ and are also division algebras at the generic point.  These Azumaya algebras have $\chi = 0$ and $c_2 = 8$. We prove that the moduli space of these Azumaya algebras is the K3 surface $X$ given by the base locus of the quadric bundle.  This is like an inverse to the classical Mukai correspondence \cite{Mukai2}.  Our original motivation for studying this problem came from mathematical physics with a view towards reformulating Mukai's results in the context of Azumaya algebras.  
%In some special cases when rank of $\Pic(X)$ and $\Pic(M)$ is 2,  we show how the Azumaya algebras trivialise \cite{IK1}. 
 
In Section 5 we extend Bridgeland's methods~\cite{Tom} to the categories of modules over Azumaya algebras and
reformulate Mukai's results on moduli spaces of sheaves in terms of derived equivalences of Azumaya algebras.  This 
also has an application to the above example in Sections 3 and 4.  Another interpretation of this equivalence is in terms of twisted sheaves \cite{Caldararu}, \cite{Huybrechts-Stellari},\cite{Yoshioka}, \cite{Huybrechts-Stellari}.
Section 6 contains some results about Chern classes of Morita equivalent Azumaya algebras. In particular we show that for an Azumaya algebra $A$ of degree $r$, if $c_2(A)$ is within $2r$ of its lower bound then $c_2(A)$ is minimal.  When
$A$ is an Azumaya algebra on a K3 surface this implies that the moduli space of these Azumaya algebras is proper \cite{Artin-deJong}.

{\it  Acknowledgements}   The second author would like to thank Bernd Kreussler for many helpful discussions.  The second author would also like to 
thank Professor Le Dung Trang for his support during the author's stay at Abdus Salam International Centre for Theoretical Physics, Trieste.
%Our original motivation for studying this problem came from mathematical physics and a view towards formulating Mukai's results in the context of Azumaya algebras.  
%We conclude with a discussion of these ideas.

\end{section}

\section{Brauer Group Basics}\label{sec:bb}
In this section we recall some basic definitions and set up the notation 
for later use.  The base field $k$ is an algebraically closed field
with characteristic p such that p is not equal to 2 or 3.  Some of the
general results also require that $p$ does not divide the index,
period, and rank of the Azumaya algebras.  We assume that all our
schemes are integral and reduced, unless specified otherwise.  We also
assume all cohomology is \'etale unless otherwise noted.  We denote by
$D(A)$ the bounded derived category of complex of $A$-modules with
coherent cohomology.
   
The Brauer group $\Br(K)$, of a field $K$ is the set of equivalence classes of central simple algebras over K. We have 
$$\Br(K) = H^2(Gal (\overbar{K}/K), \overbar{K}^*) = H^2_{\acute{e}t}(\Spec K, \G_m).$$ This notion generalises to schemes also.
For more details of this and the following results see \cite{Milne}.

%\subsection{Azumaya algebras}\label{subsec:azu}
\begin{definition} Let $X$ be a reduced integral scheme. An {\it Azumaya algebra A} on $X$ is 
a coherent sheaf of $\O_X$-algebras such that
\begin{enumerate}
\item $A$ is a locally free $\O_X$-module of rank $n^2$.
\item For every geometric point $p \in X$, the fibre
$A_p := A \otimes \kappa(p)$ is isomorphic to $\kappa(p)^{n \times n}$,
the algebra of $n \times n$ matrices over $\kappa(p)$.  
\end{enumerate}
\end{definition}
The integer $n$ is called the 
$degree$ of $A$.  The second condition is equivalent to $A_p$ being a central
simple $\kappa(p)$-algebra.  \'Etale locally an Azumaya algebra is the sheaf of endomorphisms of a vector
bundle.  

\begin{theorem}
Let $X$ be a scheme, and let $A$ be a sheaf of $\O_X$-algebras, that
is locally free of finite rank as a $O_X$-module.  Then $A$
is an Azumaya algebra on $X$ if and only if for every point 
$p \in X$, there exists an \'etale neighbourhood $\pi:U \rightarrow X$
of $p$ and a locally free sheaf $\E$ on $U$ such that 
$\pi^*A \simeq \cEnd_U(\E)$.
\end{theorem}

%\begin{definition}
%An Azumaya algebra $\A$ on a variety $X$ is a coherent sheaf 
%with a central $\O_X$-algebra structure such that the fibres 
%$\A_p:=A \otimes k(p)$ are central simple $k(p)$ algebras for all point $p$
%in $X$.
%\end{definition}

%The following result is well known, for example see \cite{AdeJ}.

\begin{proposition} Let $A$ be an Azumaya algebra.  Then
there is a trace pairing $tr : A \otimes A \rightarrow \O$
that is nondegenerate so $A \simeq A^*$ as $\O$-modules.
\end{proposition}

Two Azumaya algebras $A$ and $A'$ are said to be {\it Morita equivalent} if and only if there exist locally free sheaves 
$V, W$ such that
$$\cEnd V \otimes A  \ \simeq \  \cEnd W \otimes A'.$$  This defines an equivalence relation on the set of Azumaya algebras on a scheme X, because 
$\cEnd V \otimes \cEnd W \simeq \cEnd (V \otimes W)$.  The tensor product of two Azumaya algebras is another Azumaya algebra, and this operation is compatible with the equivalence relation.

\begin{definition}\label{def:Br(X)}
The Brauer group, $\Br(X)$ of a scheme X, is the group of Morita equivalence classes of Azumaya algebras, under the operation \\
$[A][A'] = [A \otimes A']$.  The identity element is 
$[\O_X]$ and $[A]^{-1} = [A^{\mathrm{op}}]$.
\end{definition}

\begin{definition}\label{def:period-index}
The {\it period} of an Azumaya algebra $A$, is its order in the Brauer group.  The {\it index} is $\sqrt{\dim_K D}$ where $K$ is the field of fractions of
$Z(A)$ and $A \otimes K \simeq D^{n \times n}$ with $D$ a division algebra
with centre $K$.  
\end{definition}

A central simple algebra is a division algebra if its degree equals its index.  A celebrated theorem of De Jong 
\cite{deJong} states that the index of a central simple algebra on a surface is equal to its period.  Hence to calculate the 
index of an Azumaya algebra on a surface it is enough to calculate its order in the Brauer group.
The following fact is well known.

\begin{proposition}
Let X be a regular, integral, scheme and $k(X)$ the field of rational functions on $X$. There is an injection $\Br(X) \hookrightarrow \Br(k(X))$ induced by
restricting an Azumaya algebra to the generic point. 
\end{proposition}
Since an Azumaya algebra is torsion free, there is an injective sheaf
map  $A \hookrightarrow A \otimes k(M)$.  We refer to the 
central simple algebra $A \otimes k(M)$ as $k(A)$.
The next statement follows from the Skolem-Noether Theorem.
\begin{theorem}
Let $A$ be an Azumaya algebra on $X$.  Then any automorphism of $A$ is Zariski locally an inner automorphism. 
\end{theorem}

The automorphism group of $\GL_n(\O_X)$ is 
$\PGL_n(\O_X)$, so
an Azumaya algebra $A$ corresponds to an element of $H^1(X, \PGL_n(\O_X)).$ We denote the image of $A$ in $H^1(X, \PGL_n(\O_X))$ by $c(A)$.  If $A = \cEnd V$ for some locally free sheaf $V$ then $c(A)$ is trivial in $H^1(X, \PGL_n(\O_X))$.

We have a short exact sequence
$$1 \rightarrow \O^* \rightarrow \GL_n(\O_X) \rightarrow \PGL_n(\O_X) \rightarrow 1.$$
In the associated long exact sequence, we have
$$H^1(X, \O^*) \rightarrow H^1(X, \GL_n(\O_X)) \rightarrow H^1(X, \PGL_n(\O_X))
\stackrel{\delta}{\rightarrow} H^2(X, \O^*).$$
The maps $\d$ are compatible for varying $n$, and 
$$\d(c(A) \otimes c(A')) = \d(c(A)) \otimes \d(c(A')).$$
So we get a natural injection $\Br(X) \hookrightarrow H^2(X, \O^*)$.

In many cases this is an isomorphism, in particular when $X$ is a smooth projective 
variety.

%\begin{definition}
%Let $A$ be an Azumaya algebra of degree $n$ on $X$.  The {\it Brauer-Severi variety}, $\BS(A)$, of an Azumaya algebra $A$, is the scheme parameterising rank $n$ left ideals of $A$.  Equivalently it is the scheme parameterising rank $n$ quotients of $A$ that are also left $A-$ modules.
%\end{definition}

%It is isomorphic to a $\P^{n-1}$ bundle over $X$.  The class of $\BS(A)$ in $H^1(X, \PGL_n)$ is the same as the class of $A$.

We mention some results on moduli of Azumaya algebras from \cite{Artin-deJong} which we will use later.  Let $X$ be a smooth projective surface and let $\a \in \Br(X)$  be  a Brauer class.
If $S$ is a scheme, we consider Azumaya algebras $A_S$ of degree $n$
over $X_S = X \times S$ with fixed degree $n$ and for every closed
point $s$ of $S$, the Brauer class of $A_s$, the restriction of $A_S$ to $U \times \{s\},$  is $\a$. \\

Let $\overbar{\A}$ denote the functor  where $\overbar{\A}(S)$ is the set of isomorphism classes of Azumaya algebras satisfying the conditions  above.
If $A_S$ is a family of Azumaya algebras, the Chern class $c_2(A_s)$
for $s \in S$ is a locally constant function on $S$.  Define
$\overbar{\A}_c$ to be the functor where $\overbar{\A}_c (S)$ is the
isomorphism classes of  Azumaya algebras $A_S$ over $S$ with $c_2(A_s)
= c$ for $s \in S$.  The following Theorem is due to Artin and de Jong \cite{Artin-deJong}.

\begin{theorem}\label{thm:Artin-deJong}[Theorems 8.7.6, 8.7.7]
Let  $\overbar{\A}_c$ be as above. 
\begin{enumerate}
\item If $H^1(X, \mu_n) = 0$ then the functor $\overbar{\A}_c$ is representable as an algebraic space.
\item If $c$ denotes the minimal second Chern class for the Brauer class $\a$, then the algebraic space $\overbar{\A}_c$ if non empty is proper.
\end{enumerate}
\end{theorem}

%Recall that a proper algebraic space of dimension two is a projective surface.

%\subsection{K3 surfaces}\label{subsec:k3}
We mostly work with Azumaya algebras $A$ over
an algebraic K3 surface $X$, although the results in Section 6 apply to all surfaces. 

\begin{definition}\label{defn:k3}
A K3 surface $X$ is a smooth, simply connected, complex surface with trivial
canonical bundle.
\end{definition}

 We say that $A$
is {\it simple} if $H^0(X,A) =k$.
% and 
%maybe some other restrictions that hold when $A=\cEnd V$
%is simple.  
%We 
%say that $A$ is {\it strictly simple} if $A \otimes k(X)$ 
%is a division algebra.  

\begin{proposition}
Let $X$ be a regular, integral scheme and $E$ a vector bundle on $X$. Let $A$ be an Azumaya algebra on $X$ with $r^2 = \rk A$. 
\begin{enumerate}
%\item If $A$ is strictly simple then $A$ is simple.
\item If $E$ is a simple vector bundle then $A = \cEnd E$ is simple.
%\item $\cEnd E$ is never strictly simple if $\rk E >1$.
\item If $X$ is K3 and $A$ is simple then $H^2(X,A) =k$ and $c_2(A) \geq 2r^2-2$.
\end{enumerate}
\end{proposition}
\begin{proof}
(1) This is just by definition, since $E$ is simple means $\End E = k$, so
 $h^0(A) = h^0(\cEnd E) = 1$. \\
(2) By Serre Duality we get $H^0(A) = H^2(A^{*})$.  But
$A \simeq A^{*}$ via the trace pairing 
 $A \otimes A \rightarrow \O$, so $H^2(X, A) = k$.  The Riemann-Roch Theorem implies $\chi(A) = h^0(A) - h^1(A) + h^2(A) = [\ch A.\td \T_X]_2$.
 Since $A \simeq A^{*},$ we have $c_1(A) = 0$, so $\ch A = r^2 - c_2t^2$ and 
 $\td \T_X  = 1+ 2t^2$.
 So $\chi(A) = 2r^2 -c_2 \leq h^0(A) + h^0(A^{*})$ hence $2r^2 - c_2 \leq 2.$
\end{proof}

\section{The Koszul--Clifford Algebras and K3 Surfaces}\label{sec:kc}

In this section we review a classical construction of an Azumaya algebra, calculate its invariants and prove that it has no deformations.

\subsection{Quadrics}\label{subsec:quadrics}

For the reader's convenience we review some facts about quadratic forms on vector spaces.
Let $V$ denote an $n$ dimensional complex vector space.  By a 
\emph{quadratic form}
$Q$ on $V$ we mean
a symmetric bilinear pairing
$$Q: V \times V \rightarrow k.$$
If we choose a basis $k^n \simeq V$ then  
There is an $n \times n$ symmetric matrix $S = (s_{ij})$
corresponding to $Q$ such that for any $ v,w \in V$,
$$Q(v,w) = vS{w^T}.$$

In terms of coordinates on $V$, 
$Q(v,v)$ is a homogeneous quadratic polynomial
and hence defines an associated \emph{quadric hypersurface}
$$\overbar{Q} = \{ [v] \in {\P}V : Q(v,v) = 0\}$$
in the projective space ${\P}V.$  Conversely any quadric hypersurface
$\overbar{Q}$ in ${\P}V$, defined by the zero locus of a homogeneous polynomial of degree $2$ lifts to a quadratic form $Q$ on $V$.  

We define the \emph{rank} of $Q$ to be the
rank of the matrix $S.$  If the rank of $S$ is $n$ or equivalently if $\overbar{Q}$ is smooth then $Q$ is called \emph{nondegenerate,}  otherwise $Q$ is called \emph{degenerate.}

A linear subspace $W \subset V$ is an \emph{isotropic subspace} for $Q$ if
$Q_{|_W} \equiv 0$ or equivalently ${\P}W \subset \overbar{Q}.$  If the dimension of 
$V$ is $n$ and dimension of $W$ is $\lfloor n/2 \rfloor$ then $W$ is said to be \emph{maximal isotropic.}
The \emph{isotropic Grassmannian} 
$$\IsGr(Q)=\IsGr(\lfloor n/2 \rfloor,n,Q) = \{ W \subset V : \dim W = \lfloor n/2
\rfloor, Q_{|_W} \equiv 0 \}$$
is the variety of maximal isotropic subspaces of $Q.$  Equivalently we can consider 
the quadric hypersurface
$\overbar{Q} \subset {\P}V$ in which case the the maximal isotropic subspaces correspond to families of $\P^{n-1}$ 
contained in $\overbar{Q}.$ 
If $Q$ has rank $2n$, then 
$\IsGr(n,2n,Q)$ is smooth of dimension $n(n-1)/2$ and consists of two disjoint connected components.
If $Q$ has rank $2n-1$ and dimension $2n$, then $\IsGr(n,2n,Q)$ is still smooth, of dimension 
$n(n-1)/2$ but now has only one component \cite{GH}~Chapter 6.  For example a smooth quadric surface in $\P^3$ has 
two families of lines on it whereas a quadric cone has only one family of lines.  

The Clifford algebra $\Cl(Q)$ of a quadratic form $Q$ on a vector space $V$ is defined as 
$$\Cl(Q) : = k \la V \ra / (v \otimes v = Q(v,v)), \, \text{where} \ k\la V \ra = \bigoplus V^{\otimes n}.$$
The even degree terms form a subalgebra known as the {\it even Clifford algebra} $\Cl_0(Q).$
In the next subsection we extend these notions to vector bundles.

\subsection{A Classical Construction}\label{subsec:kc-construction}

%We review a classical construction of an Azumaya algebra on $M.$  This example seems to be well known in folklore but we could not find it in the literature.  We describe it in detail for the reader's convenience. We then calculate its invariants and subsequently
%prove that it as
%no deformations.

   Let $Q_0, Q_1, Q_2$ be rank 6 quadratic
 forms on $k^{6}$ and $\overbar{Q}_0, \overbar{Q}_1, \overbar{Q}_2$ the associated quadric
 hypersurfaces in $\P^5$.  Let $\overbar{\cal{Q}}$ be the net of quadrics spanned by 
 $\overbar{Q}_0, \overbar{Q}_1$ and $\overbar{Q}_2$.  We denote by
 $X$ the base locus of $\overbar{\cal{Q}}$.  In general
 $X$ is a complete intersection and is a smooth K3 surface of degree 8 in $\P^5.$ 
We can also view $\cal{Q}$ as a quadratic form on the trivial vector bundle $\O_{\P^2}^6$, with values in $\O_{\P^2}(1)$ 
$$\cal{Q}: \Sym^2(\O_{\P^2}^6) \rightarrow \O_{\P^2}(1).$$
We denote by $\overbar{\cal{Q}}$ the vanishing locus $V(\cal{Q})$ in the corresponding projective bundle $\P^2 \times \P^5$.  
For generic choice of $Q_0,Q_1,Q_2,$ the locus of degenerate quadrics  is a smooth degree $6$ curve in $\P^2$ given by $C= V( \det (\cal{Q})),$ and furthermore
all quadrics in the net have rank $\geq 5$.  We assume that this is the case unless specified otherwise.

More generally we can define a quadratic form on a vector bundle with values in
a line bundle.

\begin{definition}\label{def:quadric-bundle}
Let $V$ be a rank $n$ vector bundle on a variety $Z$.  By a quadratic form $\cal{Q}$ on 
$V$ with values in a line bundle $\L$, we mean a map
$$\cal{Q} : \Sym^2 V \rightarrow {\L}.$$
The rank of  
$\cal{Q}$ restricted to the generic fibre is defined to be the {\it rank} of $\cal{Q}$.
The fibrewise vanishing locus $V(\cal{Q})$ in ${\P}V$, defines a family of quadric hypersurfaces over the base variety $Z$.  We call this the {\it quadric bundle} $\overbar{\cal{Q}}$ over the base $Z$. 
\end{definition}

If $\L = \O$ there is a Clifford algebra $\Cl(\cal{Q})$ associated to $(\cal{Q}, V),$

$$\Cl(\cal{Q}) = T(V) / (v \otimes v - \cal{Q}(v,v))$$ where $v$ is a  local section
of $V$ \cite{Bichsel-Knus}.  The relations defining $\Cl(Q)$ are in $\Sym^{\bullet} V$.
Locally we have 
$$\Cl(\cal{Q}) = \O \oplus V \oplus \wedge^2 V \oplus \cdots \oplus \wedge^n V.$$  The even subalgebra is
$$\Cl_0(\cal{Q}) = \O \oplus \wedge^2 V \oplus \wedge^4 V \oplus \cdots .$$

  If $\L \neq \O$, there isn't a notion of a 
Clifford algebra as such. However following the construction in \cite{Bichsel-Knus} we can still form an algebra with a sheaf structure 
$$\Cl_0(Q) = \O \oplus (\wedge^2 V \otimes \L^{-1}) \oplus (\wedge^4 V \otimes \L^{-2}) \oplus \cdots .$$

We take the tensor product $T(V, \L)$ of the $\Z$ graded algebras $ \oplus V^{\otimes i}$ and 
$\oplus \L^i$,
$$T(V, \L) := \bigoplus V^{\otimes i} \otimes_{Z} \oplus \L^{-j}.$$
$\Sym^2(V) \subset V \otimes_Z V$ so we may consider $\cal{Q}$ as a relation in $T(V, \L)$
and define the Clifford algebra $\Cl(\cal{Q})$ to be the quotient of 
$T(V, \L)$ with defining relations $\cal{Q}$.  More precisely let 
$I \subset T(V, \L)$ be the two sided ideal generated by sections of the form $t \otimes t - \cal{Q}(t),$ for $t$ a local section of $V$.  Then let $\Cl(Q) := T(V, \L)/ I$.
The algebra ${\Cl}(\cal{Q})$ is no longer bigraded but if we assign degree 1 to $V$ and degree 2 to $\L^{-1}$, 
then the relation $t \otimes t - \cal{Q}(t)$ is homogeneous of degree 2.
The degree $0$ part ${\Cl}_0(\cal{Q})$ is called the {\it even Clifford Algebra} associated to the quadratic form $\cal{Q}$ on $V$.  If $\L = \O$ we get the usual even Clifford algebra in this way.  

%\begin{definition}
%Let $\cal{Q}: \Sym^2 V \rightarrow \L$ be a quadratic form on a vector bundle of rank $n$.  let $T(V, \L)$ be the tensor algebra 
%$$T(V, \L) := \bigoplus V^{\otimes i} \otimes_{Z} \oplus \L^{-j}.$$
%Assign degree $1$ to $V$ and degree $2$ to $\L^{-1}$.
%Then the {\it even Clifford algebra} ${\Cl}_0(\cal{Q})$ of the quadric bundle $\cal{Q}$ is defined as the degree $0$ subalgebra of 
%$$T(V, \L) / (t \otimes t - \cal{Q}(t), \mbox{ for all } t \mbox{ sections of } V).$$
%\end{definition}

For our next result we need the following notions.
Let $A$ and $B$ be $\Z_2$ graded algebras.  We define a new graded algebra 
$ A \otimes B$ where the multiplication has the sign conventions 
$(a \otimes b)(a' \otimes b') = (-1)^{\deg(b)\deg(a')}aa' \otimes bb'$.  There is a natural grading on the matrix algebra $k^{2 \times 2}$ called the {\it checkerboard grading} 
$$k^{2 \times 2}_0 = \begin{bmatrix}k & 0 \\ 0 & k\end{bmatrix} \quad \quad 
k^{2 \times 2}_1 = \begin{bmatrix}0 & k \\ k & 0\end{bmatrix}.$$

The following Lemma \cite{Lam} is straightforward; 
\begin{lemma}
Let $B$ be a $\Z_2$ graded algebra.  Then there exists an isomorphism of graded algebras
$\phi:k^{2 \times 2} \otimes B \rightarrow B^{2 \times 2}$.
\end{lemma}
Let $H$ be the quadratic form given by the hyperbolic lattice
$\begin{bmatrix} 0 & 1 \\ 1 & 0 \end{bmatrix}.$  It is well known that $\Cl(H) \simeq k^{2 \times 2}$ and $\Cl_0(H) \simeq k \oplus k$. It follows that 
$\Cl(H) \otimes B \simeq B^{2 \times 2}$.  In particular $$B^{2 \times 2}_0 = \begin{bmatrix}B_0 & B_1 \\ B_1 & B_0\end{bmatrix}.$$

\begin{proposition}\label{prop:Cl_0(Q)}
Let $\cal{Q}:\Sym^2 \O_{\P^2}^6 \rightarrow \O_{\P^2}(1)$ be a quadratic form on the trivial bundle $\O_{\P^2}^6$.  
Then
 ${\Cl}_0(\cal{Q})$ is an Azumaya algebra over its centre $Z$, whose underlying scheme is 
  $\phi: M  \rightarrow \P^2$, the double cover of $\P^2$ branched along
the curve $V( \det \cal{Q}).$  
\end{proposition}
\begin{proof}
In other words
 ${\Cl}_0(\cal{Q}) = \phi_{*}A$, where $A$ is an Azumaya algebra on $M.$
This follows from local calculations of the Clifford algebra, by modifying the methods in \cite{Lam} to modules over $k\llbracket t\rrbracket$.
Since the quadratic form has rank six over any point $p$ not in $V (\det \cal{Q})$  we  
compute that 
$\Cl_0(\cal{Q}) \otimes k(p) \simeq (k^{4 \times 4})^{\oplus 2}$.
Next we need to 
compute the behaviour over a local transverse slice of the sextic $V (\det \cal{Q})$.
We may assume that 
$\cal{Q}= H \oplus H \oplus q_t$ where 
$$q_t = {\begin{bmatrix}1 & 0 \\
0 & t\end{bmatrix}}.$$
We first consider the quadratic form given by $Q_t = H \oplus q_t$.  Then% and then extend the result to $H \oplus Q_t$.
% Then by definition
%$\Cl(Q_t) = \C \la e_1, e_2, e_3, e_4 \ra / v \otimes v = Q(v,v)$, where 
%$\C \la e_1, e_2, e_3, e_4 \ra$ denotes all formal products and their linear combinations.  
$$\Cl(Q_t) = \Cl(H) \otimes \Cl(q_t) \simeq k^{2 \times 2} \otimes \Cl(q_t), \,\,\, \text{and}$$ 
$$\Cl_0(Q_t) =
\begin{bmatrix}\Cl_0(q_t) & \Cl_1(q_t) \\
\Cl_1(q_t) & \Cl_0(q_t)\end{bmatrix}.$$  Let $e_1, e_2, e_3, e_4$ denote the standard basis for the underlying $k \llbracket t\rrbracket$ module on which $Q_t$ is a quadratic form. 
Then $\Cl(q_t) = k \la e_3, e_4 \ra /( e_3^2 - 1, e_4^2 - t, e_3e_4 + e_3e_4)$ and $\Cl_0(q_t) \simeq k\llbracket\sqrt{t}\rrbracket$ since $(e_3e_4)^2=-t$.  The element $e_3$ in $\Cl_1(q_t)$ is invertible.
An easy computation shows that $e_3\!\Cl_1(q_t) = \Cl_1(q_t)e_3 = \Cl_0(q_t)$.  Moreover, conjugation by $e_3=e_3^{-1}$ acts by an automorphism on $\Cl_0(q_t)$
since $e_3 \Cl_0(q_t) e_3 = \Cl_0(q_t)$.
Now conjugation by the element 
$$\begin{bmatrix}1 & 0 \\
0 & e_3\end{bmatrix}$$ acting on $\Cl_0(Q_t)$ gives an isomorphism
$ \Cl_0(Q_t) \simeq \Cl_0(q_t)^{2 \times 2} \simeq k\llbracket\sqrt{t}\rrbracket^{2 \times 2}$.
The generalisation to $\cal{Q}$ is essentially the same argument. 

% From Lemma 3.3 it follows that
%$\Cl_0(\cal{Q}) = 
%\begin{bmatrix}\Cl_0(Q_t) & \Cl_1(Q_t) \\ \Cl_1(Q_t & \Cl_0(Q_t)\end{bmatrix}$.  
%Since $e_3$ is still an invertible element in $\Cl_1(Q_t)$, conjugation by the matrix ${\begin{bmatrix}I_2 & 0 \\ 0 & e_3I_2\end{bmatrix}}$ gives $\Cl_0(\cal{Q}) = M_4(\C\llbracket\sqrt{t}\rrbracket)$.  

So $\Cl_0(\cal{Q})$ is an Azumaya algebra over its centre $Z(A)$.  The underlying scheme
$M$  of $Z(A)$ is a double cover of $\P^2$ branched over the locus of degenerate quadrics.  
\end{proof}
%Hence
%$\Cl_0(\cal{Q})$ is an Azumaya algebra over the K3 surface $\phi: M \rightarrow \P^2$, branched over $C$. 
%$\Cl_0(\cal{Q}) \simeq {\begin{bmatrix}{\Cl_0(Q_t) & \Cl_0(Q_t) \\ \Cl_0(Q_t) & \Cl_0(Q_t)}\end{bmatrix}} \simeq M_4(\C\llbracket\sqrt{t}\rrbracket)$.

%In Section 4, we will see that
%given a K3 surface $M$ which is a double cover of a plane, branched over a smooth curve $C$, there are $2^9(2^{10}+1)$ inequivalent nets of quadrics $\cal{Q}$ such that $V(\det(\cal{Q})) = C$.  So we
%get $2^9(2^{10}+1)$ non isomorphic Azumaya algebras on $M$ which have the same invariants as $A$ above. 

%\subsection{Numerical invariants}\label{subsec:kc-invariants}

We use Koszul algebras to calculate invariants of $A$.
There are several equivalent ways of defining a Koszul algebra one of which is the following: 

\begin{definition}
Let $R=\bigoplus_{n \geq 0} R_n$ be a graded algebra with $R_0=k$.
Suppose $R$ has a minimal free resolution of $R_0=k$ with finite
generated free modules 
$$ \cdots \rightarrow R^{\oplus n_2} \stackrel{d_2}{\rightarrow} R^{\oplus n_1} 
\stackrel{d_1}{\rightarrow} R \stackrel{d}{\rightarrow} k \rightarrow 0. $$
Then $R$ is Koszul if and only if the entries of $d_i$ are all of degree one.
\end{definition}

The Koszul dual of an algebra is defined as the Yoneda algebra
$R^! := \Ext^{\bullet}(k,k).$ 

\begin{proposition}
The homogeneous coordinate ring of a complete intersection of quadrics
is Koszul.
 If $R$ is Koszul, then $H_R(t) H_{R^!}(-t) =1$ where $H_R(t)$ is the 
Hilbert series of $R$.
\end{proposition}

For a proof see \cite{Manin}.

When $R=k \la V \ra/J$ for $J \subset V \otimes V$, the Koszul dual has
the nice form $ R^! = k \la V^* \ra/ J^\perp.$
Consider the family $\cal{Q}$ of quadratic forms on $V \simeq k^{6}$ as before.  The base locus $X$ of $\overbar{\cal{Q}}$ in $\P^5$ is a complete intersection of three quadrics $\overbar{Q}_0,
\overbar{Q}_1, \overbar{Q}_2.$
Let $y = (y_0, y_1, y_2, y_3, y_4, y_5)$ denote the coordinates on $V$, 
and let $f_0, f_1, f_2$ be the respective defining polynomials of the quadric hypersurfaces.

The homogeneous coordinate ring $R$, of $X$ is Koszul and is given by

$$ R = k \la V \ra / J, \ \text{where} \, 
J := ( [y_i, y_j], {f_0}(y,y), {f_1}(y,y), {f_2}(y,y)).$$

%It is easy to write down the Hilbert series of a complete intersection $X$ in 
%$\P^n$.  Let $X=V(f_1,f_2,\ldots,f_k)$ where $f_i$ is homogeneous of 
%degree $n_i$.  Then for $R$, the homogeneous coordinate ring of $X$,
%$$ H_R(t) = \prod_i\frac{(1-t^{n_i})}{(1-t)^{n+1}}.$$
%In our case $X$ is the complete intersection of three quadrics, so
The Hilbert series of $X$ is 

$$H_R(t) =\frac{(1-t^2)^3}{(1-t)^6}=\frac{(1+t)^3}{(1-t)^3}.$$
Since $\displaystyle{H_R(t){H_{R^!}(-t)}}=1,$ it follows
that $H_R(t) = H_{R^!}(t)$.  Computing this Hilbert series yields
%We note various facts.

$$ \frac{(1+t)^3}{(1-t)^3} = -1 + \sum^\infty_{i=0} (4i^2+2)t^i$$
We also have the equality of algebras
$$ v_2(R^!) := \bigoplus_{n \geq 0} R^!_{2n} = 
\bigoplus_{n \geq 0} H^0(M, A(n)),$$ where $A$ is as in \ref{prop:Cl_0(Q)} \cite{Buchweitz} \cite{Kapranov}.
So we may interpret
$R^!$ as the polarisation of $A$ by the invertible bimodule $A(\O_M(1)/2)$,
but it will be more convenient to use the polarisation $O_M(1)$ itself.

The Hilbert Series of $v_2(R^!)$ consists of the even terms of $H_{R^!}(t)$, so 
$$H_{v_2(R^!)}(t)= -1+\sum_n (16n^2+2)t^n = 1 + 18t + 66t^2+\cdots .$$
By the Hirzebruch-Riemann-Roch formula we get $\chi(A(n))= 16n^2+32-c_2(A) = h^0(A(n))$ for $n>>0$.
%$$ \chi(A(n))  =  [\ch(\O_M(n).\ch(A).\td(\cal{T}_M)]_2 $$
%where $[-]_2$ denotes the degree two component and $\cal{T}_M$ is the tangent
%sheaf of $M$.  

%Using formulae \ref{eq:chern} and \ref{eq:td} from Section~\ref{sec:bb} 
%we get $$\displaystyle{\ch(\O_M(n))  =  1+n.ht + \frac{n^2h^2}{2} t^2 =  1+n.ht + n^2t^2.}$$  

%Since $A \simeq A^*$ via the trace pairing $A \otimes A \rightarrow \O_M$ it follows
%that
%$c_1(A) =0$ and hence $\ch(A) = 16 -c_2(A)t^2.$

%\begin{eqnarray*}
%$\chi(A(n)) & = & [(1+n.ht+n^2t^2)(16-c_2(A)t^2)(1+2t^2)]_2 \\
%& = &
%= 16n^2+32-c_2(A).$
%\end{eqnarray*}
%Since $\O_M(1)$ is ample we know that $\chi(A(n)) = h^0(A(n))$ for 
%$n >> 0$.
Comparing coefficients with the known Hilbert series $H_{v_2(R^!)}(t)$, we get $\chi(A(n)) = 16 n^2 + 2$ for all n.
%$$H_{v_2(R^!)}(t) = -1+ \sum_n (16n^2+2)t^n = 1 + 18t+66t^2+ \cdots $$
%we get that
%$\chi(A(n)) = 16n^2+2 \quad\text{for}\quad  n>>0.$
%Since $c_2(A) = 30$ 

In particular,
$\chi(A) = h^0(A) -h^1(A) + h^2(A) = 2.$
From the first term of the Hilbert series it follows that $h^0(A) = 1$,
hence $h^2(A) = 1$ and $h^1(A) = 0$. 
Recall from Section 2, that for a simple Azumaya algebra $A$, $c_2(A) \geq 2r^2 -2$. 
In our case, $\rk(A) = 16,$ and $c_2(A)=30,$ so 
this Azumaya algebra obtains its possible lower bound for $c_2(A).$

%\begin{proposition}\label{prop:deform-azu}
The first order deformations of an Azumaya algebra, $A,$ fixing its centre 
$M$, are given by $H^1(M,A/\O_M)$ with obstructions in
$H^2(M,A/\O_M)$, as described in~\cite{Ingalls}.
%\end{proposition}

%\begin{corollary}
%Let $A$ be as in Proposition \ref{prop:Cl_0(Q)}.  Then $A$ has no deformations.
%\end{corollary}
%\begin{proof}
There is a short exact sequence 
$$ 0 \rightarrow \O_M \rightarrow A \rightarrow A/\O_M \rightarrow 0$$
which is split by the trace map $tr:A \rightarrow \O_M$.
From the associated long exact sequence of cohomology it follows that
$H^1(M,A/\O_M)=H^2(M,A/\O_M)=0$. Therefore $A$
has no deformations that fix its centre.

In general $A$ is non-trivial in $\Br(M).$  This is easy to see if
$\rho(M)=1$ i.e. $\Pic(M) = \Z.h$ where $h^2=2$.  The algebra $A$ is
trivial in $\Br(M)$ if and only if $A = V \otimes V^\vee$ for some
vector bundle $V$.  Let $c_1, c_2$ be the first and second Chern
classes of $V$.  Since $\rho(M)=1$, $c_1 = n.h$ for some $n$ which
implies $c_1^2 = 2n^2$.  Then $\chi(A) = \chi(V \otimes V^\vee) = 32 +
3c_1^2 - 8c_2 = 32 + 6n^2 - 8c_2 = 2.$  This implies $n^2=4k-5$, for
some $k$.  So $n$ has to be odd.  Suppose $n=2t+1$, then simplifying
the equation yields $2t^2+2t-2k=-3$ which is a contradiction.  Later on in Section 4 we will prove using a reduction technique 
that in general $A$ has period 2.  It has period 1 if there exists a class $u \in \Pic(M)$ which has odd intersection with the hyperplane class.  This 
is related to a classical moduli problem as well \cite{Mukai2}, \cite{Ingalls-Khalid}.

We summarise these results in the following proposition.

\begin{proposition}\label{prop:deform-azu}
Let $A$ be as in Proposition \ref{prop:Cl_0(Q)}.  Then,
\begin{itemize}
\item $h^0(A) = 1, \chi(A)=2, c_2(A) = 30$,
\item $A$ has no deformations fixing its center.
\item the generic $A$ has order $2$ in $\Br(M)$.  It is trivial if and only if there exists a class $u$ in $\Pic(M)$ which has odd 
intersection with the hyperplane class.
\end{itemize}
\end{proposition}

\section{Moduli of Azumaya algebras}\label{sec:azu-moduli}

In this section starting with any K3 surface $M$, which is a double cover of $\P^2$, we construct a 
family of Azumaya algebras on $M$, parameterised by an associated K3 surface $X$ of degree $8$.  We prove that the moduli space of these Azumaya algebras is a compact, irreducible surface.  We then show that $X$ maps injectively into the moduli space hence $X$ is isomorphic to the moduli space.  These algebras have index 2 and period 2 and so are division algebras at the 
generic point.

We first describe 
%in subsection~\ref{subsec:quad-theta} 
how to associate to $M$ a net of quadrics 
$\cal{\overbar{Q}} \subset \P^5.$

%\subsection{Theta characteristics and nets of quadrics}\label{subsec:quad-theta}
\begin{definition}\label{def:theta}
A theta characteristic on a smooth curve $C$ is a line bundle $L$ such that
$L^{\otimes2} \cong \omega_{C}$.  It is called even if $H^0(L) = 
0 \ (\text{mod} \ 2)$ and 
odd otherwise. 
\end{definition}

Given an ineffective theta characteristic we get a family of symmetric
matrices $\cal{S}$ defining a quadric bundle $\cal{Q}$ with $V(\det
\cal{Q})= C$ \cite{Beauville} \cite{Turin}.
\begin{theorem}\label{thm:quad-theta}
Let $\phi: M \rightarrow \P^2$ be a double cover of $\P^2$ branched along a smooth sextic curve $C.$  Let $L$ be an even theta characteristic on $C$ such that
$H^0(L) = 0$.  
\begin{enumerate}
\item Then we have a resolution of $L$; 
$$ 0 \rightarrow {\O_{\P^2}(-2)}^6 \stackrel{\cal{S}}{\rightarrow} {\O_{\P^2}(-1)}^6 
\rightarrow {L} \rightarrow 0.$$
The matrix $\cal{S}$ is symmetric and defines
a quadric bundle $\cal{Q}$ with
$V(\det \cal{Q}) = C$. The base locus $X$ of the corresponding net $\overbar{\cal{Q}}$ in $\P^5$ is smooth.

\item Conversely, given such a matrix $\cal{S}$, its cokernel is a theta characteristic on $C$
with $H^0(L)=0$.

\end{enumerate}

\end{theorem}

 There are $2^{g-1}(2^g +1)$ inequivalent even theta characteristics on a 
curve
of genus $g$.  The curve $C$ has genus 10,  and for generic $C$ the theta
characteristics are not effective so we get  $2^9(2^{10}+1)$ different nets of quadrics for each $M$.  Since $C$ is smooth, the base locus $X$ of $\overbar{\cal{Q}}$ in
$\P^5$ is also smooth (\cite{Turin} Lemmas~1.1 and 2.6). % Recall that $X$ is a K3 surface of degree 8 in $\P^5$ realised as a
%complete intersection of three quadric hypersurfaces.  
In general $\Pic (X)$ and $\Pic (M)$ have rank 1.

If $h^0(L) = 2$,  $\Pic (M)$ has rank $2.$  Then $M$ is isomorphic
to $X$ and embeds in $\P^5$ as a degree $8$ surface.  However it is not a complete
intersection so the methods in this paper do not apply.
 We describe this case in \cite{Ingalls-Khalid}.

%\subsection{Reduction}\label{subsec:reduction}

There is a construction that allows one to pass from a quadric
bundle of dimension $d$ to a quadric bundle of dimension $d-2$,
and depends on choosing a section. 
% We apply it to our example of the quadric bundle $\cal{Q}$ of dimension $6.$  
Let $M, X$ and $\cal{Q}$ be as in Theorem~\ref{thm:quad-theta} and pick a point $x \in X.$ 
Note that $x$ is a smooth point of every quadric in the net 
(Lemma 2.6~\cite{Turin}). Let $a \in \P^2$, let $\overbar{Q}_a$ be the corresponding quadric hypersurface in $\P^5$ and
%a quadric in the net with associated matrix $S_x$  
let $T_x \overbar{Q}_a$ be the projective tangent space to $\overbar{Q}_a$ at $x.$  
Then
$\overbar{Q}_a \cap T_x \overbar{Q}_a$ is a singular quadric hypersurface in $T_x \overbar{Q}_a$
consisting of the union of all lines in $\overbar{Q}_a$ passing through $x.$

Let $F \cong \P^3$ be a
subspace of $T_x \overbar{Q}_a$ not containing $x.$
 If $a$ is not in $C$
then $\overbar{Q}_a$ is smooth, and $T_x \overbar{Q}_a \cap \overbar{Q}_a$ is a cone over
a smooth quadric surface $\overbar{Q}_x(a)$ in $F.$ If $a \in C$ then $\overbar{Q}_a$ is singular with a vertex and $T_x \overbar{Q}_a \cap \overbar{Q}_a$ is a cone over a singular quadric 
surface $\overbar{Q}_x(a)$ with a vertex.
%So given any $a$, projecting away from $x$, we get a quadric surface 
%$\overbar{Q}_x(a)$, which is singular or smooth depending on whether $a$ is in $C$ or not.
%For details of proofs see \cite{GH}~Chapter 6.
Projecting away from $x$ defines a family of quadric surfaces in a $\P^3$ bundle.  A priori it is not obvious that this 
$\P^3$ bundle is the projectivisation of a vector bundle but in fact this family of quadric surfaces is the vanishing locus of a quadratic form on the vector bundle $\Omega^1_{\P^2}(1) \oplus \O_{\P^2}^2$.

\begin{theorem}\label{thm:B_x}
Let $\phi:M \rightarrow \P^2$ be a double cover of $\P^2$ branched along a smooth sextic curve $C$ and let $X$ and $\cal{Q}$ be as in Theorem~\ref{thm:quad-theta} for some choice 
of a theta characteristic $L$ where $H^0(L) = 0.$  Let $x \in X.$ 
\begin{enumerate}
\item Then projection away from $x$ induces a quadratic form $\cal{Q}_x$ on the vector
bundle $V \simeq \Omega^1_{\P^2}(1) \oplus \O_{\P^2}^2$, with values in $\O_{\P^2}(1).$
%$$ \cal{Q}_x: \Sym^2 V \rightarrow \O_{\P^2}(1).$$  
The quadrics have rank 4 everywhere away from 
$C$ and have rank 3 on $C$.  
\item The even Clifford algebra $B_x:= \Cl_0(\cal{Q}_x)$ is a degree 2 Azumaya algebra on $M$ with $\chi(B_x) = 0$ and $c_2(B_x) = 8$.
\end{enumerate}
\end{theorem}

\begin{proof}

%Let  $(y_0,y_1, y_2, y_3, y_4, y_5)$ denote coordinates on $k^6$ and 
Let $a = [a_0:a_1:a_2]$ denote homogeneous coordinates on $\P^2$.  We may assume with out loss of generality that 
$x = [0:0:0:0:0:1]$.   Let $\cal{S}$ be the family of symmetric matrices corresponding to $\cal{Q}$ and let $T_x \cal{Q}$ be the rank $5$ subbundle of $\O_{\P^2}^6$, given by 

$$0 \rightarrow T_x{\cal{Q}} \rightarrow \O^6_{\P^2} 
\xrightarrow{\la x\cal{S}, \cdot \ra} \O_{\P^2}(1) \rightarrow 0.$$

Projection away from $x$ induces a map $\pi: \O_{\P^2}^6 \rightarrow \O_{\P^2}^5$ and we get 
a commutative diagram,

$$\begin{CD}
 % \mu_r & @=   & \mu_r  & @=   & \mu_r \\
 % @VVV &      &  @VVV &      & @VVV \\
0 @>>> T_x \cal{Q}  & @>>> & \O_{\P^2}^6     & @>{\la x \cal{S}, \cdot \ra}>> & \O_{\P^2}(1) &  @>>> & 0 \\
&      &  @VVV   &      &  @VV{\pi}V &      & @| &     & \\
0 @>>> V & @>>> & \O_{\P^2}^5 & @>{\la \pi(x \cal{S}), \cdot \ra}>> & \O_{\P^2}(1) & @>>> & 0. 
\end{CD}$$

Here $\la \pi(x \cal{S}), \cdot \ra$ denotes the induced map from $\O^5$ to $\O_{\P^2}(1)$, after projection.  After an appropriate change of basis of $\O_{\P^2}^5$ we may assume that
$\la \pi (x \cal{S}), \cdot \ra = \la (a_0, a_1, a_2, 0,0), \cdot \ra$.
The Euler exact sequence on $\P^2$ twisted with $\O_{\P^2}(1)$ is 
$$0 \rightarrow \Omega^1_{\P^2}(1) \rightarrow \O_{\P^2}^3 \stackrel{\s}{\rightarrow}
\O_{\P^2}(1) \rightarrow 0,$$
where for $\l = (\l_0, \l_1, \l_2) \in H^0(\O^3_{\P^2}),$ 
 $\s (\l) = a_0 \l_0 + a_1\l_1 + a_2\l_2.$
Then $V = \ker{\s} \oplus \O_{\P^2}^2$, i.e.
$V \simeq \Omega^1_{\P^2}(1) \oplus \O_{\P^2}^2$. 

The degenerate quadrics in $\cal Q_x$ are exactly those obtained from $\cal{Q}$ via projection away from $x$ so they are parameterised by $C$ and have rank $3$.

%\begin{theorem}\label{thm:moduli-B}
%Let $\phi: M \rightarrow \P^2$ be a K3 surface realised as a double cover of $\P^2$ branched along a smooth sextic $C$.  Let 
%$L$ be an even theta characteristic on $C$ such that $H^0(L)=0$ and $\overbar{\cal{Q}}$
%the quadric bundle obtained from a minimal resolution of $L$. Let $X$ be 
%the base locus of $\overbar{\cal{Q}}.$
%Then $X$ is a moduli space of Azumaya algebras on $M$ with
%$\chi = 0$ and rank = $4.$
%\end{theorem}

The same arguments as in Proposition~\ref{prop:Cl_0(Q)}, show that
locally across a transversal slice of the sextic we have 
$\Cl_0(\cal{Q}_x) = k \llbracket  \sqrt{t}\rrbracket^{2 \times 2}$.   So the even Clifford algebra ${\Cl}_0(\cal{Q}_x)$ is an Azumaya algebra $B_x$ on $M$, or in other words
$\phi_{*}B_x = \Cl_0(\cal{Q}_x)$.
As a sheaf $\Cl_0(\cal{Q}_x) \simeq \O \oplus (\wedge^2 V \otimes \O(-1)) \oplus (\wedge^4 V \otimes \O(-2))$.
It follows that $H^0(B_x) = k$.  A computation using the Grothendieck-Riemann-Roch formula shows that 
$\chi(B_x)= 0$ and hence $c_2(B_x) =8.$  
\end{proof}

We now show that for all $x \in X$,  $B_x$ is Morita equivalent to the Azumaya
algebra $A$ where $A$ is constructed as in 
Proposition~\ref{prop:Cl_0(Q)}.

\begin{theorem}
Let $A$ be as in Proposition \ref{prop:Cl_0(Q)} and $B_x$ be as above.  Then $A$ is Morita equivalent to $B_x$.
\end{theorem}

\begin{proof}
It is enough to show this at the generic point of $\P^2$.  Since $\P^2$ contains an affine open subset, we have
$K:=K(\P^2) = K(\C^2) = k(x,y)$, so we can view 
$\cal{Q}$ as a single quadric on $V:=k(x,y)^6$.  Pick a point $x$ in the base locus of $\cal{Q}$ and consider 
$T_x(\cal{Q}) \cap \cal{Q}$.  
This is a cone over a quadric $Q'$ of rank 4.  Since $T_x(\cal{Q}) = \{w \in V \mid \cal{Q}(w,x)=0\}$ and since $\cal{Q}(x,x)=0$ 
we get an induced quadratic form on $T_x(\cal{Q})/<x>$.  Choose a basis $f_i$ of  $T_x(\cal{Q})/<x>$ and lift to 
representatives 
$e_i$ in $T_x(\cal{Q})$.  Since $e_i \in T_x(\cal{Q})$ it follows that $\cal{Q}(e_i,x) = 0$.  It remains to find one element 
$w$ to complete $\{e_i,x\}$ to a basis of $V$.  We require that $\cal{Q}(w,e_i)= \cal{Q}(w,w) = 0$.  Since $\cal{Q}$ is non degenerate it follows that $\cal{Q}(w,x) \neq 0$.  Consider $\{s \in V \mid \cal{Q}(s,e_i) = 0\}$.  This is a two 
dimensional subspace of $V$ hence describes a line in $\P V$.  This line meets $\cal{Q}$ in two points one of which is 
$x$.  Let $w$ be the other point.  Then we have that $\cal{Q} = Q' \oplus H$ where $H$ is the hyperbolic lattice.  Therefore
by standard results on Clifford algebras \cite{Lam} we get that $\Cl_0(Q) = \Cl_0(Q')^{2x2}$ which implies 
$A_K = {B_x}^{2x2}_K$.  Since $\Br(M)$ injects into $\Br(k(M))$ it follows that $A$ is Morita equivalent to 
$B_x$.
\end{proof}
 
We now show that the Azumaya algebras $B_x$ have period $2,$ which implies they have index $2$.  %A celebrated theorem of Johan deJong 
%\cite{deJong} states that on an algebraic surface, the period of a central simple algebra is equal to its index.  It follows that they have index 2.  A central simple algebra is a division algebra if its degree equals its index.  
Hence they are division algebras at the generic point.  Since $A$ is Morita equivalent to $B_x$, it also has index 2.

For the next Theorem, we define a quadric bundle to be a variety $Q$ with a map $\pi:Q \rightarrow B$ to a smooth variety $B$.  So all fibres of $Q$ are quadrics in a fixed projective space $\P^{2n-1}=\P(V)$ and the map $\pi$ is flat and projective.  We further assume that the total space $Q$ is smooth, and that the fibres of $Q$ have rank $2n$ generically, and have rank $2n-1$ on a smooth 
divisor of the base $B$.

Let $\IsGr(Q)$ be the maximal isotropic Grassmanian of $Q$.  So 
$$\IsGr(Q) = \{ (W,b) \in \Gr(n,V) \times B \st W \subseteq Q_b \mbox{ for some } b \in B \},$$
where $\Gr(n,V)$ is the Grassmannian of $n$ dimensional subspaces of $V$.

Let $A$ be an Azumaya algebra of rank $m^2$ on a smooth variety $B$. 
The Brauer-Severi variety of $A$, $\BSV(A)$, is the variety of $m$ 
dimensional quotients of $A$ that are also left $A$-modules.  It is a $\P^{m-1}$ bundle over the 
centre of $A$ and its class in $H^1(\PGL(n))$ is the same as the class of $A$.
Note that we have a natural inclusion $$\BSV(A) \subseteq \Grass(A,m)$$
where $\Grass(A,m)$ is the Grassmanian of $m$-dimensional quotients of $A$.

\begin{theorem}\label{thm:bs(B_x)}
Let $Q\rightarrow B$ be a quadric bundle in $\P^{2n-1}=\P(V),$
then there is a natural closed immersion 
$$i:\IsGr(Q) \rightarrow \BSV(\Cl_0(Q))$$ 
$$W \mapsto (\Cl(Q)W)_0=\Cl_1(Q)W$$
which sends a maximal isotropic subspace to the ideal it generates in the even Clifford algebra.  This map is an isomorphism if and only if $n \leq 3$.
\end{theorem}
\begin{proof} 
Let $D$ be the discriminant divisor of the quadric bundle in the base $B$.
Note that the Clifford algebra of $Q$ is a vector bundle of rank $2^{2n}$ over the base $B$, away from $D$.  Hence the even Clifford algebra has rank $2^{2n-1}$.  This is not a square but the same argument as in Proposition~\ref{prop:Cl_0(Q)} showa that the even Clifford algebra, $\Cl_0(Q),$ is an Azumaya algebra of rank $2^{2n-2}$ over the 
double cover $\tilde{B}$ of $B$ ramified on $D$.
This implies that the Brauer Severi variety of $\Cl_0(Q)$ is in 
$\Gr(\Cl_0(Q),2^{n-1})$, and that $\BSV(\Cl_0(Q))\rightarrow \tilde{B}$ is a $\P^{2^{n-1}-1}$ bundle.
The dimension of the maximal isotropic Grassmannian, or spinor variety, is 
$n(n-1)/2$.  Since $n(n-1)/2=2^{n-1}-1$ only for $n=1,2,3$ this 
shows an isomorphism is impossible when $2n > 6$.
%(A possible approach would be to show that given an isotropic space $W$, the left ideal  
%$(\Cl(Q)W)_0=\Cl_1(Q)W$ has codimension $2^{n-1}$ and that different isotropic subspaces give rise to different ideals.)

A local basis computation shows that given a maximal isotropic subspace $W$,
we have that $\Cl(Q)_1W$ has codimension $2^{n-1}$ in $\Cl_1(Q)$.
To show the map is an injection, let $W,W'$ be isotropic subspaces above a point
$b$ in $B$ and restrict the quadratic form to a neighbourhood $U$ of $b$.  
We will simply write $\Cl(Q)$ for the sections of $\Cl(Q)$ over $U$.
Now suppose $(\Cl(Q)W)_0=(\Cl(Q)W')_0$.  So we can conclude that 
$\Cl(Q)_1W=\Cl(Q)_1W'$.  Hence  $\Cl(Q) \Cl(Q)_1W=\Cl(Q)\Cl(Q)_1W'$ but since
$\Cl(Q) \Cl(Q)_1=\Cl(Q),$ we see that $\Cl(Q)W=\Cl(Q)W'$.  Since we are working locally, we can filter the Clifford 
algebra in the usual way.  A basis argument using the fact that $W$ is isotropic shows that we can now 
recover $W$ from the associated graded of $\Cl(Q)W$ so we get $W=W'$.

We now show that the map is an isomorphism explicitly 
for the case $n=2$.  The case $n=3$ is similar.
This is well known if there is no discriminant, so we compute
an explicit isomorphism for the map locally at a point in the discriminant.
We take a formal disk transversal to the divisor $B$ and denote this by 
$\Delta_t = \Spec k\llbracket t \rrbracket$ with $V(t=0)=D$ locally.
We can choose a basis of $V \times \Delta$ so that our quadric has Gram matrix:
$$\begin{pmatrix} 
-1 & 0 \\
0 & t \end{pmatrix}
 \bigoplus 
\begin{pmatrix} 
0 & \frac{1}{2} \\
\frac{1}{2} & 0 \end{pmatrix}.$$

We show that we have explicit parametrizations:
$$\begin{CD} \P^1_{[a:b]} \times \Delta_{\theta} @. \P^1_{[x:y]} \times \Delta_s \\
                 @VVV   @VVV\\
                \IsGr(Q) @>{i}>> \BSV(\Cl_0(Q)) 
\end{CD}$$
that allow us to show that the map $i$ is an isomorphism.
We first parametrize $\IsGr(Q)$.

Choose a basis $e_1,\ldots,e_4$ of $V$ and consider the natural
basis $e_i \wedge e_j$ on the Plucker embedding of $\Gr(2,V)$.  Consider two arbitrary vectors in $V$ with coordinates 
$u=(u_i)$ and $v=(v_i)$, and impose 
$$-u_1^2+tu_2^2+u_3u_4=-v_1^2+tv_2^2+v_3v_4= -u_1v_1+tv_2u_2+(u_3v_4+u_4v_3)/2=0$$ to ensure that $u,v$ span an isotropic space.  This describes a closed subvariety of $\P(V) \times \P(V)$.  We now compute the image of the map $(u,v) \mapsto u \wedge v$ in $\Gr(2,V) \times \Delta$.  We write $\sum E_{ij} e_i \wedge e_j$
for a vector in $\wedge^2V$.

So we can compute that $\IsGr(Q)$ is the closed 
subvariety of $\Gr(2,V)$ described
by the equations:
$$\begin{array}{rr}
   2E_{1 3}E_{2 4}-E_{1 2}E_{3 4},&
        2E_{1 4}E_{2 3}+E_{1 2}E_{3 4},\\
        2E_{1 3}E_{1 4}+E_{3 4}^2,&
        E_{1 2}E_{1 4}+E_{2 4}E_{3 4},\\
        E_{1 2}E_{1 3}-E_{2 3}E_{3 4},&
        E_{1 2}^2-2E_{2 3}E_{2 4},\\
        tE_{2 4}^2-E_{1 4}^2,&
        2tE_{2 3}E_{2 4}-E_{3 4}^2,\\
        tE_{1 2}E_{2 4}+E_{1 4}E_{3 4},&
        tE_{2 3}^2-E_{1 3}^2,\\
        tE_{1 2}E_{2 3}-E_{1 3}E_{3 4}. 
\end{array}$$

It is easy to verify that this has the desired properties of having fibres the disjoint union of two conics for $t \neq 0$ and a single conic when $t=0$.

Let $\theta^2=t$ and let $[a:b] \in \P^1$.
Note that the vectors
$$w_1=(ab,0,a^2,b^2) \mbox{  and  }  w_2=(0,ab,-a^2\theta,b^2\theta)$$ 
span a two dimensional isotropic subspace of quadratic form $Q_t$ when 
$a \neq 0$.  

Now we can write our parametrization
$$ \P^1\times \Delta \rightarrow \IsGr(Q)$$
given by 
$$[a:b,\theta] \mapsto w_1 \wedge w_2=[ab,-a^2\theta: b^2\theta: -a^2:-b^2:2ab\theta,\theta^2]=[E_{12}:E_{13}:E_{14}:E_{23}:E_{24}:E_{34},t].$$
Note that this extends naturally to include the case $a=0$.
A further computation shows that the image of this map is also cut out by the 
same equations above proving, that we have an isomorphism.

Next we parametrize $\BSV(\Cl_0(Q))$.
Now we have another isomorphism 
$$ \Cl_0 Q \simeq k\llbracket s \rrbracket^{2 \times 2}, $$
where $s^2=t$.
Note that $\Cl_0 Q$ is generated by even products of the $e_i$ 
$$ e_1e_2 \mapsto \begin{pmatrix} s & 0 \\ 0 & -s \end{pmatrix}, \\
 e_1e_3 \mapsto \begin{pmatrix} 0 & -1 \\ 0 & 0 \end{pmatrix}, \\
 e_1e_4 \mapsto \begin{pmatrix} 0 & 0 \\ -1 & 0 \end{pmatrix}, \\
 e_2e_3 \mapsto \begin{pmatrix} 0 & 1 \\ 0 & 0 \end{pmatrix}, \\
 e_2e_4 \mapsto  \begin{pmatrix} 0 & 0 \\ s & 0 \end{pmatrix}, \\$$ $$
 e_3e_4 \mapsto  \begin{pmatrix} 1 & 0 \\ 0 & 0 \end{pmatrix}, 
e_1e_2e_3e_4  \mapsto  \begin{pmatrix} s & 0 \\ 0 & 0 \end{pmatrix}, \\
e_1e_2(e_4+e_3)(e_4-e_3)  \mapsto  \begin{pmatrix} s & 0 \\ 0 & s \end{pmatrix}. \\$$

The Brauer-Severi variety of an Azumaya algebra $A$ is the variety of left module quotients of $A$ of dimension two.  This is also naturally ideals of $A$ of codimension two.  
This shows that the Brauer-Severi variety of $\Cl_0 Q$ is isomorphic to $\P^1 \times \Delta$ with parameters $[x:y,\alpha]$ giving the left ideal of $k\llbracket s \rrbracket^{2 \times 2}$ generated by $$(xe_{11}+ye_{12},xe_{21}+ye_{22},s-\alpha)$$ where the $e_{ij}$ are the matrix units.
So we compute the image 
$$ W \mapsto (\Cl(Q) W)_0=\Cl_1(Q)W$$
where $W$ is spanned by $w_1,w_2$. 
To this end we compute $\Cl_1(Q)w_1$ and $\Cl_1(Q)w_2$ to obtain,
$$\frac{b}{a} e_4w_2 +e_2w_1-\frac{a}{b} e_3w_2+2\frac{a}{b} e_1e_2e_3w_1=e_1e_2(e_4+e_3)(e_4-e_3)+\theta \mapsto s+\theta,$$
$$e_1e_3e_4w_1 = \frac{a}{b}e_1e_3-e_3e_4 \mapsto -\frac{a}{b}e_{12}-e_{11}, $$
$$ e_4w_1= \frac{a}{b} -e_1e_4-\frac{a}{b}e_3e_4 \mapsto \frac{a}{b}e_{22}+e_{21}.$$
So we see that $w_1 \wedge w_2 \mapsto [b:a,-\theta]=[x:y,\alpha]$,
which clearly gives an isomorphism of our two parametrizations 
$$\P^1_{[a:b]} \times \Delta_\theta \rightarrow \P^1_{[x:y]} \times \Delta_s.$$
\end{proof}

The maximal isotropic Grassmannian of $\cal{Q}$ consists of two disjoint conics away from the degenerate locus.  These merge into the same conic over the degenerate locus.  It follows that 
$\BS(B_x) \simeq OGr(2,4,\cal{Q}_x)$ which is a conic bundle on $M$.  This is related to a classical correspondence 
where $M$ is a moduli space of rank 2 sheaves of $X$ with $c_1 = \O_X(1), c_2 = 4$.  These sheaves arise from 
quadrics in $\P^5$ and are known in the literature as {\it spinor bundles}.  In general $M$ is not a fine moduli space that is a universal sheaf does not exist.  However a quasi universal sheaf $\cal E$ does exist such that 
$\cal E_{|_{X \times m}} \simeq {F_m}^{2}$ where $F_m$ is a representative of the the isomorphism class of bundles corresponding to $m$ \cite{Khalid}.  

Consider the incidence correspondence in $\P^2 \times \P^5 \times Gr(2,5)$,
$$I = \{(a,x, \Lambda) \mid a \in \P^2, x \in X, \Lambda \subset T_x(\overbar{Q}_a) \cap \overbar{Q}_a \}$$

This is a conic bundle  $p: I \rightarrow X \times M$ such that for each $x \in X$, 
$I_{|_{x \times M}} \simeq OGr(2,4, \cal{Q}_x) \simeq \BS(B_x)$ and $I_{|_{X \times m}} \simeq \P F_m$.

The obstruction to the existence of a universal sheaf is
the obstruction to lifting $I$
to a rank $2$ vector bundle on $X \times M$.  In fact $I$ admits
a section if and only if there exists an element in $\Pic(X)$ which has odd intersection with $\O_X(1)$.   In particular if 
$\rho(M) =1$ then $I$ does not admit a section therefore the period of $B_x$ is $2$.  if $X$ contains a line then $X \simeq M$ and $I$ lifts to 
a vector bundle on $X \times M$.  In this case the moduli problem is symmetric \cite{Ingalls-Khalid}.  

Consider the short exact sequence as in \cite{Chan-Ingalls}

$$0 \rightarrow  \O \rightarrow J \rightarrow T_{I/X \times M} \rightarrow 0,$$
where $T_{I/X \times M}$ is the relative tangent sheaf.
Then $\B = \pi_*(\cEnd J) \simeq \pi_*J$ is an Azumaya algebra on $X \times M$ and 
$\BS(\B) \simeq I$.  So $\B$ is a universal family on $X \times M$ such that $\B_{|_{x \times M}} \simeq B_x$.   Moreover 
$\B_{|_{X \times m}} \simeq {F_m} \otimes {F_m}^{\vee}$, up to tensoring by line bundles, where $F_m$ is a spinor bundle on $X$ corresponding to $m \in M$.
%Curiously, $\B$ is also isomorphic to the above mentioned quasi-universal sheaf 
%$\cal E,$ up to tensoring by line bundles from the base.  
As an application of Theorem~5.8 we get a derived equivalence 
$D(X) \equiv D(B)$ where $B = \pi_{M_*} \cal End\,{\cal{E}}$.  The
underlying gerbe of $B$ is the same as that of $B_x$ as defined in~\ref{underGerbe}.

We use the results from 
\cite{Chan-Ingalls} to compute the Chern invariants of $\BS(B_x)$.

\begin{eqnarray*}
{\chi_{top}}(\BS(B_x)) & = & 2 \chi_{top}(M) = 48, \\
\chi(\O_{\BS(B_x)}) & = & \chi(\O_M) = 2, \\
K_{\BS(B_x)}^3 & = & c_2(B_x)  = 8  
\end{eqnarray*}

%chi(T_{\BS(B_x)}) = -16$.
%The original motivation for studying this conic bundle came from an example due to 
%Mukai~\cite{Mukai2}.  
%In general $M$ is a non fine moduli space of sheaves on $X$ because the conic bundle $I$ does not lift to a vector bundle.  In some special cases $M$ is a fine moduli space and the moduli problem is actually symmetric \cite{IK1}.
%Theorem~\ref{thm:moduli-B} is like an inverse of this moduli problem. However there is a subtlety here because we have to make a choice of a theta characteristic.  So whereas given an $X$ in $\P^5$ there is a unique $M$, given an $M$ there isn't a unique corresponding K3 surface $X$. There are $2^9(2^{10}+1)$ even theta characteristics so there are $2^9(2^{10}+1)$ choices for $X$ all of which are non isomorphic as K3 surfaces.  The obstruction to the existence of a universal sheaf is an element $\alpha \in Br(M)$ which can be described in terms of a Hodge isometric embedding $T_X \rightarrow T_M$ via the Fourier-Mukai map on cohomology.
%There are $2^9(2^{10}+1)$ even sublattices of $T_M$, of index $2$, which arise as transcendental lattices of a K3 surface $X$ \cite{BertVanGeeman}.  So given a K3 surface $M$ we get 
%$2^9(2^{10}+1)$ non isomorphic division algebras on $M$ whose moduli spaces are the corresponding distinct K3 surfaces $X$.

We need one more result to proceed.  Given a curve $S$ of genus $2$, there is a pencil of quadrics $\cal{X}$ in $\P^5$ such that $S$ is isomorphic to the Stein factorisation of the family of planes in $\cal{X}$.  The base locus $U$ of $\cal{X}$ is a complete intersection of two quadrics.  Given any point $u \in U,$ and a quadric $Q$ in $\cal{X}$, the variety
$T_u(Q) \cap Q$ is a cone over a quadric surface.  So projecting away from $u$ we get a quadric surface bundle on $S$.  The variety of lines in this quadric surface bundle is a conic bundle $P$ on $S$.  Varying $u$ we get a conic bundle $\cal{P}$ on $U \times S$.
 %which corresponds to an element in the Brauer group $\Br(U \times S) = H^2(U \times S, \O^{*})$.  As Newstead points out in~\cite{Newstead}, this Brauer group is trivial .
 Since $U \times S$ has torsion free cohomology and $H^2(U\times S,\O)=0$ this conic bundle lifts to 
 a universal family $\cal F$ of rank 2 vector bundles of odd degree on $S$ \cite{Newstead}.   We summarise the construction and some results of~\cite{Newstead} in the
following Theorem.

\begin{theorem}
Let $S$ be a genus two curve.  Then the moduli space of rank 2 bundles with fixed determinant of odd degree on $S$  can be constructed as the base locus $U$ of a pencil of quadrics $\cal{X}$ in $\P^5$.  The universal rank 2 vector bundle of the moduli space gives a conic bundle over $U \times S$.  For a given $u$ in $U$, this conic bundle over $\{ u \} \times S$ is the variety of lines in the projection from $u$ of the singular quadric $T_u(Q) \cap Q$.  In particular, given $p,q$ in $U$ the corresponding conic bundles are isomorphic if and only if $p=q$.
\end{theorem}

%\begin{theorem}\label{thm:morita}
%SKIP THIS FOR NOW BECAUSE WE DON'T HAVE PROOF
%Let $M$, $\cal{Q}$,  $X$ be as in Theorem~\ref{thm:quad-theta} and $\cal{A}$ the Azumaya algebra
%on $M$ given by 
%$\Cl_0{\cal{Q}} = \phi_{*}{\cal{A}}.$
%Let $\B_p$ be an Azumaya algebra on $M$ as in Proposition~\ref{prop:B_p}.
%Then  $\B_p$ is Morita equivalent to $\cal{A}.$
%\end{theorem}
%\begin{proof}

%\end{proof}

The above construction can also be studied using Koszul duality and Clifford algebras, but we will not undertake this point of view in this paper. 

\begin{theorem}\label{thm:moduli-B}
Let $\phi: M \rightarrow \P^2$ be a K3 surface realised as a double cover of $\P^2$ branched along a smooth sextic $C$.  Let 
$L$ be an even theta characteristic on $C$ such that $H^0(L)=0$ and
let $\bar{\cal{Q}}$
the net of quadrics in $\P^5$ obtained from a minimal resolution of $L$. Let $X$ be 
the base locus of $\bar{\cal{Q}}.$
Then $X$ is a moduli space of rank $4$ 
Azumaya algebras on $M$ with $\chi = 0$, and $c_2 = 8$. 
\end{theorem}

\begin{proof}
For each $x \in X$, we construct a rank $4$ Azumaya algebra $B_x$ with
$\chi(B_x) = 0$ and $c_2(B_x) = 8$, as in Theorem~\ref{thm:B_x}.   Let
$\M$ denote the moduli stack of rank 4 Azumaya algebras on $M$ with
$c_2 = 8$ and same underlying gerbe as $B_x$ for some $x \in X$ as
defined in~\ref{underGerbe}. 
Since $h^0(B_x) = 1$, it follows that $h^1(B_x) = 2$.  From the long exact cohomology sequence associated to the short exact sequence 
$$0 \rightarrow \O \rightarrow B_x \rightarrow B_x/\O \rightarrow 1 $$
we get $h^1(B_x/\O) = 2$ and $h^2(B_x/\O) = 0$.  Therefore by
Proposition~\ref{prop:deform-azu}, we have that $\M$ is two dimensional.

%In addition for all $p$, $B_p$ has degree 2 and index 2 so
%$B_p \otimes k(M)$ is a division algebra.  It seems natural then that $X$ should be a component of the moduli stack of Azumaya algebras on $M$ with same underlying gerbe as $B_p$.  However we need to check some details.
 
Recall that an Azumaya algebra $A$ is simple if $h^0(A) = 1$.
Given a simple Azumaya algebra $A$ of rank $r^2$ (equivalently degree
$r$) on a K3 surface we know that $c_2(A) \geq 2r^2 -2$.  So for
simple rank 4 Azumaya algebras we get $c_2 \geq 6$.
 By Corollary~\ref{cor:min-c_2} it follows that the Azumaya algebras $B_x$ have minimal $c_2$, so $\M$ is a smooth projective surface (see Theorem~\ref{thm:Artin-deJong}).

The underlying gerbe of $B_x$ is an element of the
discrete group
$H^2(M, \mu_2),$ as defined in~\ref{underGerbe} hence stays constant as we vary $x$ in $X$. 
So we have a morphism 
$\psi: X \rightarrow \M$.  We need to show that it is injective to conclude that $X \simeq \M$.  It is enough to show that
for distinct points $x,y$ in $X$, $\BS(B_x)$ is not isomorphic to $\BS(B_y)$.

By Theorem~\ref{thm:bs(B_x)}, we have that $\BS(B_x),$ the
Brauer-Severi variety of $B_x$, is isomorphic to the maximal isotropic
Grassmannian of $\IsGr(2,4, \cal{Q}_x)$.  Let $S$ be a smooth genus 2
curve in $M,$ which is an element of the linear system $\O_M(1)$.
Then $S$ corresponds to a pencil of quadrics $\cal{X}$ contained in
$\cal{Q}$. The base locus $U$ of this pencil is a moduli space of rank
2 vector bundles of odd degree on $S$ as in Theorem~4.6.  Each point
$x \in U$ defines a conic bundle $P_x$ on $S$.  When $x \in X \subset
U$, we have that
$P_x \simeq \IsGr(2,4, \cal{Q}_{x_{|_S}})$ and so, by Theorem~\ref{thm:bs(B_x)}, is isomorphic to $\BS(B_x)_{|_S}$.

Suppose  there are two distinct points $x$ and $y$ in $X$ such that the Brauer-Severi varieties 
$\BS(B_x)$ and $\BS(B_y)$ are isomorphic.   Then
$\BS(B_x)_{|_S} \simeq \BS(B_y)_{|_S}$ are isomorphic conic bundles which contradicts Theorem 4.6.  Therefore the map $\psi: X \rightarrow \M$ is injective and $X \simeq \M$.

\end{proof}

In general $M$ is a non fine moduli space of sheaves on $X$ because the conic bundle $I$ does not lift to a vector bundle.  
Theorem~\ref{thm:moduli-B} forms a type of  inverse of this moduli problem. However there is a subtlety here because we have to make a choice of a theta characteristic.  So whereas given an $X$ in $\P^5$ there is a unique $M$, given an $M$ there isn't a unique corresponding K3 surface $X$. There are $2^9(2^{10}+1)$ even theta characteristics so there are $2^9(2^{10}+1)$ choices for $X$ all of which are non isomorphic as K3 surfaces.  The obstruction to the existence of a universal sheaf is an element $\alpha \in Br(M)$ which can be described in terms of a Hodge isometric embedding $T_X \rightarrow T_M$ via the Fourier-Mukai map on cohomology.
There are $2^9(2^{10}+1)$ even sublattices of $T_M$, of index $2$, which arise as transcendental lattices of a K3 surface $X$ \cite{BertVanGeeman}.  So given a K3 surface $M$ we get 
$2^9(2^{10}+1)$ non isomorphic division algebras on $M$ whose moduli spaces are the corresponding distinct K3 surfaces $X$.

An equivalent formulation of the inverse of the Mukai correspondence
is in terms of twisted sheaves, so $X$ can be viewed as a moduli space
of twisted sheaves on $M$ as in \cite{Caldararu}, \cite{Huybrechts-Stellari} \cite{Yoshioka}.

\section{Derived Equivalence}

Let $X$ be a K3 surface.  For the product of $X \times M$
we let $\pi_X, \pi_M$ be the projections, and we denote the fibres over points $m \in M$ and $x \in X$ by $X_m, M_x$.  Let $v$ be a Mukai vector such that the moduli space $M$ of vector bundles with Mukai
vector $v$ is a surface.  The following results due to Mukai are
well known \cite{Mukai1}, \cite{Mukai2}.  

\begin{proposition}
The moduli space $M$ of stable vector bundles 
$E_m$ with Mukai vector $v$, if non empty and compact, is a K3 surface. There is a quasi-universal
sheaf $\E$ on $X \times M$ such that for each point $m$ of 
$M$ we have $\E \otimes \O_{X_m} \simeq E_m^{\oplus n}$.
\end{proposition} 

From the quasi-universal sheaf $\E$ we construct an Azumaya algebra on $M$.
\begin{proposition}
Let $B = {\pi_M}_* \cEnd_{X \times M} \E$.  Then $B$
is an Azumaya algebra on $M$ and ${\pi_M}_* \E$ is a
$B$-module and $\E$ is a $\pi_M^* B$ module.
\end{proposition}
\begin{proof}
Since each $E_m$ is a simple, we know that $\End(E_m) =k$,
and so 
$$H^0(X_m,\cEnd( \E) \otimes \O_{X_m}) \simeq 
H^0(X_m,\E \otimes \E^* \otimes \O_{X_m}) $$ $$\simeq 
H^0(X,E_m^{\oplus n} \otimes (E_m^{\oplus n})^\vee) \simeq 
\End(E_m^{\oplus n},E_m^{\oplus n}) \simeq k^{n \times n}$$
and since $B$ is a coherent sheaf of $\O_X$-algebras, by cohomology and 
base change we see that $B$ is Azumaya.
\end{proof}

Let $A$ be an $\O_X$-algebra and let $B$ be an $\O_M$-algebra.
By this we mean that $\O_X$ is central in $A$.  Let $B^o$
be the opposite algebra of $B$.
As usual we define $A \boxtimes B^o = \pi_X^* A \otimes \pi_M^* B^o$.
Following Van den Bergh, we define an $A-B$ bimodule to be
a quasi-coherent sheaf $\E$ on $X \times M$ that has the 
structure of an $A \boxtimes B^o$ module.  We sometimes
write ${}_A\E_B$ to emphasise the bimodule structure. We define 
functors
$$ \E \otimes_B - = \pi_{X*}( \E \otimes_{\pi_M^*B^o} \pi^*_M -)$$
from the category of left $B$-modules to right $A$-modules and 

$$ - \otimes_A \E = \pi_{M*}( \pi^*_X- \otimes_{\pi_X^*A} \E)$$
from the category of right $A$-modules to left  $B$-modules.

We use the methods of Bridgeland \cite{Tom} to determine conditions
for these functors to induce a derived equivalence.
Let  
$X$ and $M$ be K3 surfaces, $A = \O_X$ and $B = \pi_{M*} \cEnd \E$ as above.

The main fact needed to generalise Bridgeland's proof is the following
easy Lemma.

\begin{lemma}
There is a natural isomorphism
$$ k^{1 \times n} \otimes_{k^{n \times n}} k^{n \times 1} \simeq  k.$$
\end{lemma}
\begin{proof}
We can verify this directly or use the Yoneda Lemma and the adjunction 
of tensor and $\Hom.$
\end{proof}

For every point $m$ of $M$ the fibre of $B$ is a matrix algebra
$B_m \simeq k^{n \times n}$.  There is a simple $B$-module 
$V_m \simeq k^{n \times 1}$
supported at $m$ that is unique up to isomorphism.  Note
that the modules $V_m$ are not the fibres of a vector bundle
in any way that is consistent with their structures of a $B$-module
unless $B$ is trivial in the Brauer Group $\Br M$.
Let $\Omega = \{ V_m \suchthat m \in M \}$. 
Following \cite{Tom} we call $\Omega$ a spanning class for $D(B)$ if
for all $i,m$ we have 
$ \Hom_{D(B)}(V_m, C[i]) =0 $ then $C =0$ and 
if for all $i,m$ we have $\Hom_{D(B)}(C[i],V_m)=0 $
then $C \simeq 0$.  
Since every coherent $B$-module has a non-trivial surjection
to some simple module $V_m$, we have the following Lemma.
%repeat the proof in \cite{Tom}.

\begin{lemma}
The modules $\Omega$ form a spanning class for the derived category
$D(B)$.
\end{lemma}
\begin{proof}
First note that since $\w_B \simeq B$ the two conditions are
equivalent.  So we focus our attention on the latter.
We know that for any coherent $B$-module there is some $m$
so that $\Hom_{B}(F,V_m) \neq 0,$ since $F$ must have
a simple quotient.  So if $\Hom_B(F,V_m) = 0$ for all $m$
then $F =0$.  Now let $C$ be an object of $D(B)$ such that $C \neq 0$.
Then there is a spectral sequence 
$$ \Ext^p_B(\H^q(C),V_m) \Rightarrow \Hom_{D(B)}(C[q-p],V_m).$$
Since we are in the bounded derived category there is a maximal
$q$ such that $\H^q(C) \neq 0$, so for some $m \in M$ 
there is a homomorphism in $\Hom_B(H^q(C),V_m) \neq 0$.  
Since we have $p=0$ and $q$ is maximal we have that $q-p$ is maximal,
and so $\Hom_B(H^q(C),V_m)$ can not vanish in the spectral sequence. 
\end{proof}

We need another Lemma concerning an adjunction that we prove
using the Serre functor.  Recall that if we let $f^*: \O_M \rightarrow
B$ be the structure map then the dualising bimodule for $B$ is
$$ \w_B := f^! \w_X := \Hom_B(B, \w_X ) = B^\vee \otimes \w_X$$
which satisfies Serre duality
$$ \Ext^i_B(a,b)^* = \Ext^{n-i}_B(b,\w_B \otimes a).$$
In particular $ S_B := \omega_B[n] \otimes_B -$ defines the
Serre functor on $D(B)$ so we have the equation
$$ \Hom_{D(B)}(a,b)^* = \Hom_{D(B)}(b,S(a)).$$
In our particular case where $M$ is K3 and $B$ is Azumaya,
we have that $\w_B \simeq B$
and $\w_X \simeq \O_X$ as bimodules.
% but we will not 
%use this fact until later.

\begin{lemma}  For an $\pi_M^*B$-module $a$ and a $B$-module $b$,
we have the adjunction
$$ \Hom_{D(B)}(R\pi_{M*} a,b) = 
\Hom_{D(\pi^*_M}(a,\omega_X[n] \otimes \pi^*_M b).$$
\end{lemma}
\begin{proof}
We start with the left hand side and apply the Serre functor to 
derive that 
$$\Hom_{D(B)}(R\pi_{M*} a,b) 
\simeq \Hom_{D(B)}(b,\omega_B[n] \otimes_B R\pi_{M*} a)^*$$
$$\simeq \Hom_{D(B)}(\omega_B^\vee[-n] \otimes_B b,R\pi_{M*} a)^*
\simeq \Hom_{D(\pi^*_M B)} (\pi^*_M(\omega_B^\vee[-n] \otimes_B b),a)^*$$
$$ \simeq \Hom_{D(\pi^*_M B)} (\pi^*_M\omega_B^\vee[-n] \otimes\pi^*_Mb,a)^*
$$ $$ \simeq \Hom_{D(\pi^*M B)}(a,\w_X[n] \boxtimes \w_B[n] \otimes  (\pi^*_M\omega_B^\vee[-n] \otimes\pi^*_Mb) $$
$$ \simeq \Hom_{D(\pi^*M B)}(a,\w_X[n] \otimes \pi^*_Mb).$$
\end{proof}

Let $F: D(B) \rightarrow D(X)$ be the derived functor of $\E \otimes_B -$.
Note that despite the notation, this is a right derived functor. 
Since $\E$ is a flat $B$-module we only need to derive the $\pi_{X*}$.
So 
$$ DF(-) = R\pi_{X*}(\E \otimes{\pi^*B^o} -).$$

By results of Bridgeland, we only need verify that the following conditions are satisfied.
 
\begin{theorem}
The derived functor of $F$ induces an equivalence of derived categories
if
\begin{enumerate}
\item $$ \Hom_X(FV_{m},FV_{m}) = k. $$
\item $$ \Ext^i_X(E_{m_1},E_{m_2}) =0 \mbox{ for all } i \mbox{ and } m_1 \neq m_2$$.
\item $$ E_m \otimes \omega_X \simeq E_m.$$
\end{enumerate}
\end{theorem}
\begin{proof}
We calculate  
$$FV_m = R\pi_{X*} ( \E \otimes_{\pi_M^*B} \pi_M^*V_m)$$
$$ \simeq  R\pi_{X*} ( E_m^{\oplus n} \otimes_{\O_X^{n \times n}} \pi_M^*k^{n \times 1})
$$ $$ \simeq R\pi_{X*} E_m \simeq E_m.$$  Since each $E_m$ is simple it follows that
$\Hom_X(E_m,E_m) = k$.
Since $M$ is a moduli space of sheaves on $X$, (2) is satisfied.
Condition (3) is a consequence of the fact that $X$ is a K3 surface and hence $\omega_X \simeq \O_X$.

\end{proof}

\begin{lemma}
Let 
$$FL = R(\Hom(\pi_M^*B, \E) \otimes \pi_X^* \omega_X [\dim X] \otimes_X -).$$
and
$$FR = R(\Hom(\pi_M^*B, \E) \otimes \pi_M^* \omega_B [\dim M] \otimes_X -).$$
then $FR$ and $FL$ are left and right adjoints to $DF = R(\E \otimes_B -)$.
\end{lemma}

So by a straightforward adaptation of the methods in \cite{Tom} we have the following theorem.

\begin{theorem}
Let $M$ be a K3 surface which is a moduli space of stable sheaves on a K3 surface $X$, with Mukai vector $v$.  Let $\E$ be a quasi-universal sheaf on $X \times M$. Let 
$B = \pi_{M*}(\cEnd \E)$.  Let $D(B)$ denote the derived category of left $B-modules$ on $M$.  Then there is a  derived equivalence $D(B) \simeq D(X)$.
\end{theorem}

%We apply this to our earlier example where
%$X$ is a degree 8 K3 surface in $\P^5$ given as the base locus of a net of quadrics $\cal{Q}$ and $M$ a double cover of $\P^2$ branched over $V(\det(\cal{Q}))$.  A quasi universal sheaf $\E$
%exists for this problem such that $\E_{|_{X_m}} \simeq E_m^2$ \cite{Khalid}.  Then 
%$B:= \pi_{M*}(\cEnd \E)$ is a rank $4$ Azumaya algebra on $M$, and the functor 
%$DF(-) = R\pi_{X*}(\E \otimes{\pi^*B^o} -)$ induces a derived equivalence $D(B) \simeq D(X)$.

For an equivalent formulation in terms of twisted sheaves see \cite{Caldararu}, \cite{Huybrechts-Stellari}.

%Suppose we know the following theorem.
%\begin{conjecture}
%Let $\A$ be an Azumaya algebra on a K3 surface $X$.  Fix numerical
%invariants of $\A$-modules such that there exists a moduli
%space $M$ of stable locally projective $\A$-modules with these invariants.
%Then 
%\begin{enumerate}
%\item The moduli space $M$ of locally projective $\A$-modules $E_m$ exists.
%\item There exists a quasi-universal sheaf $\E$ of $\A$-modules 
%so that $\E \otimes \O_{X_m} \simeq E_m^n$.
%\item The moduli space $M$ is a K3 surface.
%\end{enumerate}
%\end{conjecture}

%Then we obtain the following corollary.
%\begin{corollary}
%Let $B = \pi_{M*} \End_{\pi_X^*A} \E.$  Then $B$ is an Azumaya 
%algebra and $\E$ induces a derived equivalence of $D(A) \simeq D(B).$
%\end{corollary}
%The proof will essentially be the same as above.

\section{Moduli Spaces}

Let $A$ be a simple Azumaya algebra of rank $r^2$ over $X$.  Consider $M,$ 
the moduli stack of Azumaya algebras that are Morita equivalent to $A$.
Let $M_A$ be a component of $M$ that includes $A$.  The gerbe underlying $A$
in $H^2(X,\mu_r)$ is constant over $M_A$ since this group is discrete.

In this Section we show that if $A$ and $A'$ are Azumaya algebras of degree $r$ that 
have the same underlying gerbe, then $2r \mid c_2(A)-c_2(A')$.  This shows
that if the second Chern class jumps down in $M_A$ then it must jump down by at
least $2r$.  Hence if $A$ is simple and has its second Chern class within 
$2r$ of a lower bound, then it has minimal $c_2$ and by Theorem~\ref{thm:Artin-deJong}, $M_A$ is a proper moduli space.  

%Let $M$ be a K3 surface, and let $A$ be a simple Azumaya algebra over $M$ of 
%rank $r^2$.  From now on we always assume our Azumaya algebras are simple.

%\begin{proposition}
%Let $\A$ be a simple Azumaya algebra of rank $r^2$ over a K3 surface.
%\begin{enumerate}
%\item $h^0(\A) = h^2(\A) =1$
%\item $\A \simeq \A^*$
%\item $c_1(\A) =0$
%\item $\chi(\A) = 2r^2-c_2(A)$
%\item $c_2(\A) \geq 2r^2-2.$
%\item $\O$ is a direct summand of $\A$.
%\end{enumerate}
%\end{proposition}
%\begin{proof}
%The trace pairing shows that $\A \simeq \A^*$, and consequently $c_1(A)=0$.
%The first statement is the definition of an Azumaya algebra being simple,
% with statement 2 and Serre duality on a K3 surface.
%Statement 4 follows from a Rieman-Roch calculation and statement 5 
%comes from statements 4 and 1.  The last result follows from the fact
%that the trace map splits the natural inclusion of $\O \rightarrow A$.
%\end{proof}

%We now demand that $\chi(A) =0$
%or equivalently $c_2(\A) =2r^2$ or $h^1(\A) =2$.
%For such algebra $A$ we know $h^0(A/\O) = h^2(A/\O) =0$
%and $h^1(A/\O)=2$ so deformation theory tells us that
%their moduli space forms a smooth surface.  

%If we also know
%that $A$ has minimal second chern class then we know that 
%the moduli space $M_A$ will be proper and therefore projective
%since it is a surface.

We need the following diagrams:

\begin{lemma}
There is a commutative diagram
$$\begin{CD} 
\G_m & @=   & \G_m  &      &          \\
@AAA &      & @AAA  &      &          \\
\G_m & @>>> & \GL_r & @>>> & \PGL_r   \\
@AAA &      & @AAA  &      &  @AA{~}A \\
\mu_r & @>>> & \SL_r & @>>> & \PSL_r 
\end{CD}$$
of short exact sequences.
There is also a commutative diagram
$$\begin{CD}
  \mu_r & @=   & \mu_r  & @=   & \mu_r \\
  @VVV &      &  @VVV &      & @VVV \\
\SL_r  & @<<< & G     & @>>> & \G_m \\
@VVV   &      &  @VVV &      & @VV{(-)^r}V \\
\PGL_r & @<<< & \GL_r & @>{\det^{-1}}>> & \G_m
\end{CD}$$
where the bottom squares are both Cartesian.
\end{lemma}

The map $(-)^r : \G_m \rightarrow \G_m$ is the $r^{th}$ power map.
The proof is a straight forward diagram chase.  So the group $G$
above can defined as either of
$$ G = \GL \times_{\PGL} \SL \simeq \GL \times_{\G_m} \G_m.$$

So we have a commutative diagram 
$$\begin{CD}
 H^1(\PSL_r) & @>>> & H^2(\mu_r) \\
@| &           & @VVV \\
H^1(\PGL_r) & @>>>   & H^2(\G_m) 
\end{CD}$$

\begin{definition}\label{underGerbe}
Given an Azumaya algebra we obtain a class in $H^1(M,\PGL_r)$
which we can map to an element of $H^2(\mu_r)$.  We call
this the {\it underlying gerbe} of the Azumaya algebra.  The image in
$H^2(M,\G_m)$ is the corresponding element in the Brauer group.
\end{definition}

The previous Lemma shows that given a curve $C$ we have a commutative diagram.
$$ \begin{CD}
H^1(C,\GL_r) & @>{\det^{-1}}>> & \Pic C \\
@VVV & & @VVV \\
H^1(C,\PGL_r) & @>>> & H^2(C,\mu_r) 
\end{CD}$$

So if we have a vector bundle $V$ on $C$ with rank $r$ and degree $v$ then the value of 
$-v \pmod{r}$ is an invariant of $V \otimes V^{*}$ and can be considered to be the image of 
the composition
$$H^1(C,\PGL_r) \stackrel{\simeq}{\rightarrow}
 H^1(C,\PSL_r) \rightarrow H^2(C,\mu_r).$$

To proceed we need several results from \cite{Artin-deJong}.

\begin{definition}
Let $A$ be an Azumaya algebra on a smooth surface $X$.  Let $Y$ be a smooth
curve in $X$, the restriction $A \otimes \O_Y$ splits by Tsen's Theorem
so there is a vector bundle $V$ on $Y$ such that 
$A \otimes \O_Y\ \simeq \cEnd(V) \simeq V \otimes V^*$.  
Let $F \subset V$ be a subsheaf of $V$ which is locally a direct summand.
There are canonical surjective maps of right $\A$-modules:
$$ A \rightarrow A_Y \simeq V^* \otimes V \rightarrow F^* \otimes V.$$
Let $K_A$  be the kernel of the composed map, so we have an exact sequence
of right $A$-modules.
$$ 0 \rightarrow K \rightarrow A \rightarrow F^* \otimes V \rightarrow 0.$$
Artin and de Jong, \cite{Artin-deJong} 8.2, define $A' = A'(A,Y,V,F) := \cEnd K_A$ to be the 
{\it elementary transformation} of $A$.
\end{definition}

\begin{proposition}{\cite{Artin-deJong} 8.2.11}
Let $A,A'$ be two Azumaya algebras in the same central simple algebra $k(A)$
over a smooth projective surface.  Then there is a smooth curve $Y \subset X$
and a local direct summand $F \subset V$ such that $A'$ is the elementary
transformation $A'=A'(A,Y,V,F)$.  Furthermore $F$ can be taken to be an 
invertible sheaf.  Furthermore we may also replace $Y$ 
by a smooth curve linearly equivalent to a a multiple of $Y$.
\end{proposition}

Let $Q=V/F,$ be the quotient bundle on $Y$,
 and let us denote the ranks and degrees
of $V,F,Q$ as $v_0,v_1,f_0,f_1,q_0,q_1,$ respectively.  Note that
$v_0 = f_0+q_0 = r$ and $v_1 = f_1+q_1$.  We need the following formula:

\begin{proposition}{\cite{Artin-deJong} 8.3.1}
Let $A'=A'(A,Y,V,F)$ be an elementary transformation of $A$,
then
$$   c_2(A') -c_2(A)=  -f_0q_0C^2 +2(f_0q_1-f_1q_0) $$
\end{proposition}
Note: The formula as stated in the manuscript \cite{Artin-deJong} has a minor error.

\begin{proposition}
Let $A$ and $A'$ be Morita equivalent Azumaya algebras 
of degree $r$ with the same underlying gerbe in $H^2(X,\mu_r)$.  
Then $2r \mid c_2(A)-c_2(A')$.
\end{proposition}
\begin{proof}
By Proposition~7.3 \cite{Artin-deJong}, $A'$ is the elementary transformation of $A$
over a smooth curve $C$ in $M$.  This implies that there is the data
of a decomposition $A|_C = V \otimes V^*$ and a subvector bundle 
$F$ of $V$ of rank one,
from which one can build $A'$.  
We let $Q$ = $V/F$ and let $q_0,q_1,f_0,f_1,r,v_1$ be 
ranks and degrees of $Q, F$ and $V$ respectively.  So $q_0 +f_0 =r$ and $q_1 + f_1=v_1$.
We have that $f_0=1$ and $q_0=r-1$.
Note that the class of $C$ in $H^2(X,\mu_r)$ is the difference
of the classes of $A,A'$ in $H^2(X,\mu_r)$ by \cite{Artin-deJong} Corollary 8.2.10.
We can also choose
$C$ sufficiently large, as long its Chern class is zero in 
$H^2(X,\mu_r)$.  So we choose $C=rD$ such that $2r \mid C^2$.
Then (see Proposition~6.4)
\begin{eqnarray*}  
\lefteqn{  c_2(A') -c_2(A)} \\
& = & -f_0q_0C^2 +2(f_0q_1-f_1q_0) \\
                   & = & -f_0q_0C^2 + 2(v_1-f_1r) \\
 & = &  -(r-1)C^2 +2(v_1-f_1r).
\end{eqnarray*}
%  Since $\chi(\A) = 0$, it follows that $c_2(\A) = 2r^2$.  But also recall that the minimum possible $c_2(\A') = 2r^2-2$.
%so we must have 
%$$ 0< -(r-1)C^2 +2(v_1-f_1r) \geq -2. WHY DO YOU HAVE 0<$$
%Since we have that $2r \mid C^2$ we may divide by two to get
%$$ 0< -(r-1)C^2/2 + (v_1-f_1r)\geq -1$$.
So if we can show that $r \mid v_1$ then we are done. 

The value of
$-v_1 \pmod{r}$ is an invariant of $V \otimes V^*$ 
over $C$ and is the image
of the composition
$$H^1(C,\PGL_r) \stackrel{\simeq}{\rightarrow}
 H^1(C,\PSL_r) \rightarrow H^2(C,\mu_r).$$
Recall that our Azumaya algebra $A$ gives us an element of 
$$H^1(X,\PGL_r) \stackrel{\simeq}{\rightarrow} H^1(X,\PSL_r)
\rightarrow H^2(X,\mu_r) \rightarrow H^2(C,\mu_r).$$
The image of $A$ in $H^2(C, \mu_r)$ via the above mapping
$H^2(X,\mu_r) \rightarrow H^2(C,\mu_r)$, takes the value $-v_1 \pmod{r}$.

Now since $C = rD$, we conclude
that the image of $H^2(X,\mu_r) \rightarrow H^2(C,\mu_r)$
is $0$.  So we must have that $r \mid v_1$.  
\end{proof}

\begin{corollary}\label{cor:min-c_2}
Let $A$ be an Azumaya algebra of degree $r$ and suppose
that there is a lower bound $c_2(A') \geq k$ for all
Morita equivalent Azumaya algebras $A'$ of degree $r$.
If $c_2(A) < k + 2r$  then $A$ has minimal $c_2$.
\end{corollary}

\begin{proof}
Let $A'$ be an Azumaya algebra of degree $r$, Morita equivalent to $A$.  Suppose there exists an $A'$ such 
that $c_2(A') < c_2(A)$.  
Then by the above Proposition $2r \mid c_2(A) - c_2(A')$ which implies that 
$c_2(A) \geq k + 2r$ which is a contradiction.  Therefore $A$ must have minimal $c_2$. 
\end{proof}

\end{document}